\documentclass[12pt]{article}

\usepackage{amsmath,amssymb,amsthm,amscd,a4wide,upref}

\usepackage[enableskew]{youngtab}

\usepackage{hyperref}
\hypersetup{colorlinks,citecolor=blue,filecolor=black,linkcolor=blue,urlcolor=blue}
\usepackage{enumerate}
\usepackage{mathtools}
\usepackage{enumitem}

\usepackage{lmodern}     
\usepackage[T1]{fontenc}

\begin{document}


\newcommand{\ad}{{\rm ad}}
\newcommand{\cri}{{\rm cri}}
\newcommand{\row}{{\rm row}}
\newcommand{\col}{{\rm col}}
\newcommand{\Ann}{{\rm{Ann}\ts}}
\newcommand{\End}{{\rm{End}\ts}}
\newcommand{\Rep}{{\rm{Rep}\ts}}
\newcommand{\Hom}{{\rm{Hom}}}
\newcommand{\Mat}{{\rm{Mat}}}
\newcommand{\ch}{{\rm{ch}\ts}}
\newcommand{\chara}{{\rm{char}\ts}}
\newcommand{\diag}{{\rm diag}}
\newcommand{\st}{{\rm st}}
\newcommand{\non}{\nonumber}
\newcommand{\wt}{\widetilde}
\newcommand{\wh}{\widehat}
\newcommand{\ol}{\overline}
\newcommand{\ot}{\otimes}
\newcommand{\la}{\lambda}
\newcommand{\La}{\Lambda}
\newcommand{\De}{\Delta}
\newcommand{\al}{\alpha}
\newcommand{\be}{\beta}
\newcommand{\ga}{\gamma}
\newcommand{\Ga}{\Gamma}
\newcommand{\ep}{\epsilon}
\newcommand{\ka}{\kappa}
\newcommand{\vk}{\varkappa}
\newcommand{\si}{\sigma}
\newcommand{\vs}{\varsigma}
\newcommand{\vp}{\varphi}
\newcommand{\ta}{\theta}
\newcommand{\de}{\delta}
\newcommand{\ze}{\zeta}
\newcommand{\om}{\omega}
\newcommand{\Om}{\Omega}
\newcommand{\ee}{\epsilon^{}}
\newcommand{\su}{s^{}}
\newcommand{\hra}{\hookrightarrow}
\newcommand{\ve}{\varepsilon}
\newcommand{\pr}{^{\tss\prime}}
\newcommand{\ts}{\,}
\newcommand{\vac}{\mathbf{1}}
\newcommand{\vacu}{|0\rangle}
\newcommand{\di}{\partial}
\newcommand{\qin}{q^{-1}}
\newcommand{\tss}{\hspace{1pt}}
\newcommand{\Sr}{ {\rm S}}
\newcommand{\U}{ {\rm U}}
\newcommand{\cK}{ {\check{K}}}
\newcommand{\BL}{ {\overline L}}
\newcommand{\BE}{ {\overline E}}
\newcommand{\BP}{ {\overline P}}
\newcommand{\BQ}{ {\overline Q}}
\newcommand{\BR}{ {\overline R}}
\newcommand{\BT}{ {\overline T}}
\newcommand{\AAb}{\mathbb{A}\tss}
\newcommand{\CC}{\mathbb{C}\tss}
\newcommand{\KK}{\mathbb{K}\tss}
\newcommand{\QQ}{\mathbb{Q}\tss}
\newcommand{\SSb}{\mathbb{S}\tss}
\newcommand{\TT}{\mathbb{T}\tss}
\newcommand{\ZZ}{\mathbb{Z}\tss}
\newcommand{\DY}{ {\rm DY}}
\newcommand{\X}{ {\rm X}}
\newcommand{\Y}{ {\rm Y}}
\newcommand{\Z}{{\rm Z}}
\newcommand{\ZX}{{\rm ZX}}
\newcommand{\ZY}{{\rm ZY}}
\newcommand{\Ac}{\mathcal{A}}
\newcommand{\Lc}{\mathcal{L}}
\newcommand{\Mc}{\mathcal{M}}
\newcommand{\Pc}{\mathcal{P}}
\newcommand{\Qc}{\mathcal{Q}}
\newcommand{\Rc}{\mathcal{R}}
\newcommand{\Sc}{\mathcal{S}}
\newcommand{\Tc}{\mathcal{T}}
\newcommand{\Bc}{\mathcal{B}}
\newcommand{\Ec}{\mathcal{E}}
\newcommand{\Fc}{\mathcal{F}}
\newcommand{\Gc}{\mathcal{G}}
\newcommand{\Hc}{\mathcal{H}}
\newcommand{\Uc}{\mathcal{U}}
\newcommand{\Vc}{\mathcal{V}}
\newcommand{\Wc}{\mathcal{W}}
\newcommand{\Yc}{\mathcal{Y}}
\newcommand{\Cl}{\mathcal{C}l}
\newcommand{\Ar}{{\rm A}}
\newcommand{\Br}{{\rm B}}
\newcommand{\Ir}{{\rm I}}
\newcommand{\Fr}{{\rm F}}
\newcommand{\Jr}{{\rm J}}
\newcommand{\Or}{{\rm O}}
\newcommand{\GL}{{\rm GL}}
\newcommand{\Spr}{{\rm Sp}}
\newcommand{\Rr}{{\rm R}}
\newcommand{\Zr}{{\rm Z}}
\newcommand{\gl}{\mathfrak{gl}}
\newcommand{\middd}{{\rm mid}}
\newcommand{\ev}{{\rm ev}}
\newcommand{\Pf}{{\rm Pf}}
\newcommand{\Norm}{{\rm Norm\tss}}
\newcommand{\oa}{\mathfrak{o}}
\newcommand{\spa}{\mathfrak{sp}}
\newcommand{\osp}{\mathfrak{osp}}
\newcommand{\f}{\mathfrak{f}}
\newcommand{\se}{\mathfrak{s}}
\newcommand{\g}{\mathfrak{g}}
\newcommand{\h}{\mathfrak h}
\newcommand{\n}{\mathfrak n}
\newcommand{\m}{\mathfrak m}
\newcommand{\z}{\mathfrak{z}}
\newcommand{\Zgot}{\mathfrak{Z}}
\newcommand{\p}{\mathfrak{p}}
\newcommand{\sll}{\mathfrak{sl}}
\newcommand{\agot}{\mathfrak{a}}
\newcommand{\bgot}{\mathfrak{b}}
\newcommand{\qdet}{ {\rm qdet}\ts}
\newcommand{\Ber}{ {\rm Ber}\ts}
\newcommand{\HC}{ {\mathcal HC}}
\newcommand{\cdet}{{\rm cdet}}
\newcommand{\rdet}{{\rm rdet}}
\newcommand{\tr}{ {\rm tr}}
\newcommand{\gr}{ {\rm gr}\ts}
\newcommand{\str}{ {\rm str}}
\newcommand{\loc}{{\rm loc}}
\newcommand{\Gr}{{\rm G}}
\newcommand{\sgn}{ {\rm sgn}\ts}
\newcommand{\sign}{{\rm sgn}}
\newcommand{\ba}{\bar{a}}
\newcommand{\bb}{\bar{b}}
\newcommand{\bc}{\bar{c}}
\newcommand{\eb}{\bar{e}}
\newcommand{\fb}{\bar{f}}
\newcommand{\hba}{\bar{h}}
\newcommand{\bi}{\bar{\imath}}
\newcommand{\bj}{\bar{\jmath}}
\newcommand{\bk}{\bar{k}}
\newcommand{\bl}{\bar{l}}
\newcommand{\bell}{\bar{\ell}}
\newcommand{\bp}{\bar{p}}
\newcommand{\hb}{\mathbf{h}}
\newcommand{\Sym}{\mathfrak S}
\newcommand{\fand}{\quad\text{and}\quad}
\newcommand{\Fand}{\qquad\text{and}\qquad}
\newcommand{\For}{\qquad\text{or}\qquad}
\newcommand{\for}{\quad\text{or}\quad}
\newcommand{\grpr}{{\rm gr}^{\tss\prime}\ts}
\newcommand{\degpr}{{\rm deg}^{\tss\prime}\tss}
\newcommand{\bideg}{{\rm bideg}\ts}

\renewcommand{\theequation}{\arabic{section}.\arabic{equation}}

\numberwithin{equation}{section}

\newtheorem{thm}{Theorem}[section]
\newtheorem{lem}[thm]{Lemma}
\newtheorem{prop}[thm]{Proposition}
\newtheorem{cor}[thm]{Corollary}
\newtheorem{conj}[thm]{Conjecture}
\newtheorem*{mthm}{Main Theorem}
\newtheorem*{mthma}{Theorem A}
\newtheorem*{mthmb}{Theorem B}
\newtheorem*{mthmc}{Theorem C}
\newtheorem*{mthmd}{Theorem D}

\theoremstyle{definition}
\newtheorem{defin}[thm]{Definition}

\theoremstyle{remark}
\newtheorem{remark}[thm]{Remark}
\newtheorem{example}[thm]{Example}
\newtheorem{examples}[thm]{Examples}

\newcommand{\bth}{\begin{thm}}
\renewcommand{\eth}{\end{thm}}
\newcommand{\bpr}{\begin{prop}}
\newcommand{\epr}{\end{prop}}
\newcommand{\ble}{\begin{lem}}
\newcommand{\ele}{\end{lem}}
\newcommand{\bco}{\begin{cor}}
\newcommand{\eco}{\end{cor}}
\newcommand{\bde}{\begin{defin}}
\newcommand{\ede}{\end{defin}}
\newcommand{\bex}{\begin{example}}
\newcommand{\eex}{\end{example}}
\newcommand{\bes}{\begin{examples}}
\newcommand{\ees}{\end{examples}}
\newcommand{\bre}{\begin{remark}}
\newcommand{\ere}{\end{remark}}
\newcommand{\bcj}{\begin{conj}}
\newcommand{\ecj}{\end{conj}}

\newcommand{\bal}{\begin{aligned}}
\newcommand{\eal}{\end{aligned}}
\newcommand{\beq}{\begin{equation}}
\newcommand{\eeq}{\end{equation}}
\newcommand{\ben}{\begin{equation*}}
\newcommand{\een}{\end{equation*}}

\newcommand{\bpf}{\begin{proof}}
\newcommand{\epf}{\end{proof}}

\def\beql#1{\begin{equation}\label{#1}}

\newcommand{\Res}{\mathop{\mathrm{Res}}}

\title{\Large\bf Gaussian generators for the Yangian associated with\\
the Lie superalgebra $\osp(1|2m)$}

\author{Alexander Molev\ \   and\ \   Eric Ragoucy}

\date{\it To the memory of Georgia Benkart}

\maketitle


\begin{abstract}
We give a new presentation of the Yangian for the orthosymplectic
Lie superalgebra $\osp_{1|2m}$. It relies on the
Gauss decomposition of the generator matrix
in the $R$-matrix presentation. The defining relations between the Gaussian
generators are derived from a new version of the Drinfeld-type presentation of
the Yangian for $\osp_{1|2}$ and some additional relations
in the Yangian for $\osp_{1|4}$ by an application of the embedding theorem
for the super-Yangians.



%

\end{abstract}


%

\section{Introduction}\label{sec:int}
\setcounter{equation}{0}

The {\em Yangian} $\Y(\osp_{N|2m})$ associated with the orthosymplectic Lie superalgebra
$\osp_{N|2m}$ is a deformation of the universal enveloping algebra
$\U(\osp_{N|2m}\tss[u])$ in the class of Hopf algebras. The original definition
in terms of an $R$-matrix presentation and basic
properties of the Yangian are due to Arnaudon {\it et al.\/}~\cite{aacfr:rp}.
Drinfeld-type presentations of the Yangian $\Y(\osp_{N|2m})$
and extended Yangian $\X(\osp_{N|2m})$ with $N\geqslant 3$
were constructed in a recent work \cite{m:dt}. Our goal in this paper is to produce similar
presentations in the case $N=1$ (Theorem~\ref{thm:dp} and Corollary~\ref{cor:modpy}).

It is well-known that the Yangians associated with simple Lie algebras admit a few
presentations which are suitable for different applications in representation theory
and mathematical physics. In particular, the {\em Drinfeld presentation} originated in
\cite{d:nr} was essential
for the classification of the finite-dimensional irreducible representations.

Explicit isomorphisms between the $R$-matrix and Drinfeld presentations
of the Yangians associated with the classical Lie algebras were produced in
\cite{bk:pp} and \cite{jlm:ib}. In the case of the super Yangian
for the general linear Lie superalgebra, such an isomorphism between
the $R$-matrix presentation of \cite{n:qb} and a Drinfeld-type presentation of \cite{s:yl}
was given in \cite{g:gd}; see also \cite{p:pp} and \cite{t:sa} for generalizations
to arbitrary Borel
subalgebras.

A key role in the above-mentioned constructions is played by the Gauss decomposition
of the generator matrix of the (super) Yangian, which yields a presentation in terms of
the {\em Gaussian generators}. We use the same approach for the Yangians
associated with $\osp_{1|2m}$
in this paper, and
our arguments rely on the {\em embedding theorem} proved in \cite{m:dt}. It allows one
to regard the Yangian $\Y(\osp_{1|2m-2})$ as a subalgebra of $\Y(\osp_{1|2m})$, and the same
holds for their extended versions.
Therefore, a significant part of calculations is reduced to those in the algebras $\Y(\osp_{1|2})$
and $\Y(\osp_{1|4})$.

A Drinfeld-type presentation of the Yangian $\Y(\osp_{1|2})$ was given in \cite{aacfr:sy}
with the use of certain Serre-type relations. We give a different
version of this presentation involving some additional generators, but avoiding Serre-type
relations (Theorem~\ref{thm:odp} and Corollary~\ref{cor:odpy}).

The finite-dimensional irreducible representations of the algebras $\X(\osp_{1|2m})$ and $\Y(\osp_{1|2m})$
were classified in \cite{m:ry}. We apply our results to derive the classification theorem
in terms of the new presentation of the Yangian $\Y(\osp_{1|2m})$ (Proposition~\ref{prop:fdhw}).

After we posted the first version of the paper in the arXiv, we were informed by Alexander Tsymbaliuk
of the work \cite{ft:pp}, where similar results will be presented as a part of a more general project
involving presentations of the orthosymplectic Yangians associated with arbitrary parity
sequences. In particular, some closely related versions of Theorems~\ref{thm:odp} and \ref{thm:dp}
are proved in \cite{ft:pp}. We are grateful to Alexander for the illuminating discussion of those results and
their connection with the work \cite{aacfr:sy}; see Corollary~\ref{cor:serre} below.

\section{Definitions and preliminaries}
\label{sec:def}
\setcounter{equation}{0}

Introduce the
involution $i\mapsto i\pr=2m-i+2$ on
the set $\{1,2,\dots,2m+1\}$.
Consider the $\ZZ_2$-graded vector space $\CC^{1|2m}$ over $\CC$ with the
canonical basis
$e_1,e_2,\dots,e_{2m+1}$, where
the vector $e_i$ has the parity
$\bi\mod 2$ and
\ben
\bi=\begin{cases} 1\qquad\text{for}\quad i=1,\dots,m,m',\dots,1',\\
0\qquad\text{for}\quad i=m+1.
\end{cases}
\een
The endomorphism algebra $\End\CC^{1|2m}$ is then equipped with a $\ZZ_2$-gradation with
the parity of the matrix unit $e_{ij}$ found by
$\bi+\bj\mod 2$. We will identify
the algebra of
even matrices over a superalgebra $\Ac$ with the tensor product algebra
$\End\CC^{1|2m}\ot\Ac$, so that a square matrix $A=[a_{ij}]$ of size $2m+1$
is regarded as the element
\beql{matpro}
A=\sum_{i,j=1}^{2m+1}e_{ij}\ot a_{ij}(-1)^{\bi\tss\bj+\bj}\in \End\CC^{1|2m}\ot\Ac,
\eeq
where the entries $a_{ij}$ are assumed to be homogeneous of parity $\bi+\bj\mod 2$.
The involutive matrix {\em super-transposition} $t$ is defined by
$(A^t)_{ij}=a_{j'i'}(-1)^{\bi\bj+\bj}\tss\ta_i\ta_j$,
where we set
\ben
\ta_i=\begin{cases} \phantom{-}1\qquad\text{for}\quad i=1,\dots,m+1,\\
-1\qquad\text{for}\quad i=m+2,\dots,2m+1.
\end{cases}
\een
This super-transposition is associated with the bilinear form on the space $\CC^{1|2m}$
defined by the anti-diagonal matrix $G=[g_{ij}]$ with $g_{ij}=\de_{ij'}\tss\ta_i$.

A standard basis of the general linear Lie superalgebra $\gl_{1|2m}$ is formed by elements $E_{ij}$
of the parity $\bi+\bj\mod 2$ for $1\leqslant i,j\leqslant 2m+1$ with the commutation relations
\ben
[E_{ij},E_{kl}]
=\de_{kj}\ts E_{i\tss l}-\de_{i\tss l}\ts E_{kj}(-1)^{(\bi+\bj)(\bk+\bl)}.
\een
We will regard the orthosymplectic Lie superalgebra $\osp_{1|2m}$
associated with the bilinear form defined by $G$ as the subalgebra
of $\gl_{1|2m}$ spanned by the elements
\ben
F_{ij}=E_{ij}-E_{j'i'}(-1)^{\bi\tss\bj+\bi}\ts\ta_i\ta_j.
\een

Introduce the permutation operator $P$ by
\ben
P=\sum_{i,j=1}^{2m+1} e_{ij}\ot e_{ji}(-1)^{\bj}\in \End\CC^{1|2m}\ot\End\CC^{1|2m}
\een
and set
\ben
Q=\sum_{i,j=1}^{2m+1} e_{ij}\ot e_{i'j'}(-1)^{\bi\bj}\ts\ta_i\ta_j
\in \End\CC^{1|2m}\ot\End\CC^{1|2m}.
\een
The $R$-{\em matrix} associated with $\osp_{1|2m}$ is the
rational function in $u$ given by
\ben
R(u)=1-\frac{P}{u}+\frac{Q}{u-\ka},\qquad \ka=-m-\frac{1}{2}.
\een
This is a super-version of the $R$-matrix
originally found in \cite{zz:rf}.
Following \cite{aacfr:rp}, we
define the {\it extended Yangian\/}
$\X(\osp_{1|2m})$
as a $\ZZ_2$-graded algebra with generators
$t_{ij}^{(r)}$ of parity $\bi+\bj\mod 2$, where $1\leqslant i,j\leqslant 2m+1$ and $r=1,2,\dots$,
satisfying defining relations \eqref{RTT} below.
Introduce the formal series
\beql{tiju}
t_{ij}(u)=\de_{ij}+\sum_{r=1}^{\infty}t_{ij}^{(r)}\ts u^{-r}
\in\X(\osp_{1|2m})[[u^{-1}]]
\eeq
and combine them into the square matrix $T(u)=[t_{ij}(u)]$; cf. \eqref{matpro}.
Consider the elements of the tensor product algebra
$\End\CC^{1|2m}\ot\End\CC^{1|2m}\ot \X(\osp_{1|2m})[[u^{-1}]]$ given by
\ben
T_1(u)=\sum_{i,j=1}^{2m+1} e_{ij}\ot 1\ot t_{ij}(u)(-1)^{\bi\tss\bj+\bj}\fand
T_2(u)=\sum_{i,j=1}^{2m+1} 1\ot e_{ij}\ot t_{ij}(u)(-1)^{\bi\tss\bj+\bj}.
\een
The defining relations for the algebra $\X(\osp_{1|2m})$ take
the form of the $RTT$-{\em relation}
\beql{RTT}
R(u-v)\ts T_1(u)\ts T_2(v)=T_2(v)\ts T_1(u)\ts R(u-v).
\eeq

As shown in \cite{aacfr:rp}, the products $T(u-\ka)\ts T^{\tss t}(u)$
and $T^{\tss t}(u)\ts T(u-\ka)$ are scalar matrices with
\beql{ttra}
T(u-\ka)\ts T^{\tss t}(u)=T^{\tss t}(u)\ts T(u-\ka)=c(u)\tss 1,
\eeq
where $c(u)$ is a series in $u^{-1}$.
All its coefficients belong to
the center $\ZX(\osp_{1|2m})$ of $\X(\osp_{1|2m})$ and freely generate the center;
this can be derived
analogously to
the Lie algebra case considered in \cite{amr:rp}.

The {\em Yangian} $\Y(\osp_{1|2m})$
is defined as the subalgebra of
$\X(\osp_{1|2m})$ which
consists of the elements stable under
the automorphisms
\beql{muf}
t_{ij}(u)\mapsto \vp(u)\ts t_{ij}(u)
\eeq
for all series
$\vp(u)\in 1+u^{-1}\CC[[u^{-1}]]$.
As in the non-super case \cite{amr:rp}, we have the tensor product decomposition
\beql{tensordecom}
\X(\osp_{1|2m})=\ZX(\osp_{1|2m})\ot \Y(\osp_{1|2m});
\eeq
see also \cite{gk:yo}.
The Yangian $\Y(\osp_{1|2m})$ is isomorphic to the quotient
of $\X(\osp_{1|2m})$
by the relation $c(u)=1$.

An explicit form of the defining relations \eqref{RTT} can be written
in terms of the series \eqref{tiju} as follows:
\begin{align}
\big[\tss t_{ij}(u),t_{kl}(v)\big]&=\frac{1}{u-v}
\big(t_{kj}(u)\ts t_{il}(v)-t_{kj}(v)\ts t_{il}(u)\big)
(-1)^{\bi\tss\bj+\bi\tss\bk+\bj\tss\bk}
\non\\
{}&-\frac{1}{u-v-\ka}
\Big(\de_{k i\pr}\sum_{p=1}^{2m+1}\ts t_{pj}(u)\ts t_{p'l}(v)
(-1)^{\bi+\bi\tss\bj+\bj\tss\bp}\ts\ta_i\ta_p
\label{defrel}\\
&\qquad\qquad\quad
{}-\de_{l j\pr}\sum_{p=1}^{2m+1}\ts t_{k\tss p'}(v)\ts t_{ip}(u)
(-1)^{\bi\tss\bk+\bj\tss\bk+\bi\tss\bp}\ts\ta_{j'}\ta_{p'}\Big).
\non
\end{align}
In this formula and in what follows,
square brackets denote super-commutator
\ben
[a,b]=ab-ba\tss(-1)^{p(a)p(b)}
\een
for homogeneous elements $a$ and $b$ of parities $p(a)$ and $p(b)$.

The assignments
\ben
t_{ij}(u)\mapsto t_{ij}(u+c)\quad\text{with}\quad c\in \CC,
\Fand
t_{ij}(u)\mapsto t_{ji}(-u)(-1)^{\bi\tss\bj+\bj}
\een
define automorphisms of $\X(\osp_{1|2m})$ \cite{aacfr:rp}. We will need
their composition (with $c=1$) which defines
another automorphism
\beql{sigma}
\si: t_{ij}(u)\mapsto t_{ji}(-u-1)(-1)^{\bi\tss\bj+\bj}.
\eeq
The assignment
\beql{tauanti}
\tau: t_{ij}(u)\mapsto t_{ji}(u)(-1)^{\bi\tss\bj+\bj}
\eeq
defines an anti-automorphism. The latter property
is understood in the sense that
\ben
\tau(ab)=\tau(b)\tau(a)(-1)^{p(a)p(b)}
\een
for homogeneous elements $a$ and $b$ of the Yangian. Note that the maps
$\si$ and $\tau$ are not involutive but each of $\si^4$ and $\tau^4$ is the identity map.

The universal enveloping algebra $\U(\osp_{1|2m})$ can be regarded as a subalgebra of
$\X(\osp_{1|2m})$ via the embedding
\beql{emb}
F_{ij}\mapsto \frac12\big(t_{ij}^{(1)}-t_{j'i'}^{(1)}(-1)^{\bj+\bi\bj}\ts\ta_i\ta_j\big)(-1)^{\bi}.
\eeq
This fact relies on the Poincar\'e--Birkhoff--Witt theorem for the orthosymplectic Yangian
which was pointed out in \cite{aacfr:rp} and a detailed proof is given in \cite{gk:yo}; cf.
\cite[Sec.~3]{amr:rp}.
It states that the associated graded algebra
for $\Y(\osp_{1|2m})$ is isomorphic to $\U(\osp_{1|2m}[u])$.
The algebra $\X(\osp_{1|2m})$ is generated by
the coefficients of the series $c(u)$ and $t_{ij}(u)$ with the conditions
\ben
\bal
i+j&\leqslant 2m+2\qquad \text{for}\quad i=1,\dots,m,m',\dots,1'\fand\\
i+j&< 2m+2\qquad \text{for}\quad i=m+1.
\eal
\een
Moreover, given any total ordering
on the set of these generators, the ordered monomials with the powers of odd generators
not exceeding $1$, form a basis of the algebra.

The extended Yangian $\X(\osp_{1|2m})$ is a Hopf algebra with the coproduct
defined by
\beql{Delta}
\De: t_{ij}(u)\mapsto \sum_{k=1}^{2m+1} t_{ik}(u)\ot t_{kj}(u).
\eeq
For the image of the series $c(u)$ we have $\De:c(u)\mapsto c(u)\ot c(u)$ and so the Yangian
$\Y(\osp_{1|2m})$ inherits the Hopf algebra structure from $\X(\osp_{1|2m})$.

\section{Gaussian generators}
\label{sec:gd}

Let $A=[a_{ij}]$ be a $p\times p$ matrix over a ring with $1$.
Denote by $A^{ij}$ the matrix obtained from $A$
by deleting the $i$-th row
and $j$-th column. Suppose that the matrix
$A^{ij}$ is invertible.
The $ij$-{\em th quasideterminant of} $A$
is defined by the formula
\ben
|A|_{ij}=a_{ij}-r^{\tss j}_i(A^{ij})^{-1}\ts c^{\tss i}_j,
\een
where $r^{\tss j}_i$ is the row matrix obtained from the $i$-th
row of $A$ by deleting the element $a_{ij}$, and $c^{\tss i}_j$
is the column matrix obtained from the $j$-th
column of $A$ by deleting the element $a_{ij}$; see
\cite{gr:dm}.
The quasideterminant $|A|_{ij}$ is also denoted
by boxing the entry $a_{ij}$,
\ben
|A|_{ij}=\begin{vmatrix}a_{11}&\dots&a_{1j}&\dots&a_{1p}\\
                                   &\dots&      &\dots&      \\
                             a_{i1}&\dots&\boxed{a_{ij}}&\dots&a_{ip}\\
                                   &\dots&      &\dots&      \\
                             a_{p1}&\dots&a_{pj}&\dots&a_{pp}
                \end{vmatrix}.
\een

Apply the Gauss decomposition
to the generator matrix $T(u)$ associated with the extended Yangian $\X(\osp_{1|2m})$:
\beql{gd}
T(u)=F(u)\ts H(u)\ts E(u),
\eeq
where $F(u)$, $H(u)$ and $E(u)$ are uniquely determined matrices of the form
\ben
F(u)=\begin{bmatrix}
1&0&\dots&0\ts\\
f_{21}(u)&1&\dots&0\\
\vdots&\vdots&\ddots&\vdots\\
f_{1'1}(u)&f_{1'2}(u)&\dots&1
\end{bmatrix},
\qquad
E(u)=\begin{bmatrix}
\ts1&e_{12}(u)&\dots&e_{11'}(u)\ts\\
\ts0&1&\dots&e_{21'}(u)\\
\vdots&\vdots&\ddots&\vdots\\
0&0&\dots&1
\end{bmatrix},
\een
and $H(u)=\diag\ts\big[h_1(u),\dots,h_{1'}(u)\big]$.
The entries
of the matrices $F(u)$, $H(u)$ and $E(u)$ are given by well-known formulas
in terms of quasideterminants \cite{gr:tn};
see also \cite[Sec.~1.11]{m:yc}. We have
\beql{hmqua}
h_i(u)=\begin{vmatrix} t_{1\tss 1}(u)&\dots&t_{1\ts i-1}(u)&t_{1\tss i}(u)\\
                          \vdots&\ddots&\vdots&\vdots\\
                         t_{i-1\ts 1}(u)&\dots&t_{i-1\ts i-1}(u)&t_{i-1\ts i}(u)\\
                         t_{i\tss 1}(u)&\dots&t_{i\ts i-1}(u)&\boxed{t_{i\tss i}(u)}\\
           \end{vmatrix},\qquad i=1,\dots,1',
\eeq
whereas
\beql{eijmlqua}
e_{ij}(u)=h_i(u)^{-1}\ts\begin{vmatrix} t_{1\tss 1}(u)&\dots&t_{1\ts i-1}(u)&t_{1\ts j}(u)\\
                          \vdots&\ddots&\vdots&\vdots\\
                         t_{i-1\ts 1}(u)&\dots&t_{i-1\ts i-1}(u)&t_{i-1\ts j}(u)\\
                         t_{i\tss 1}(u)&\dots&t_{i\ts i-1}(u)&\boxed{t_{i\tss j}(u)}\\
           \end{vmatrix}
\eeq
and
\beql{fijlmqua}
f_{ji}(u)=\begin{vmatrix} t_{1\tss 1}(u)&\dots&t_{1\ts i-1}(u)&t_{1\tss i}(u)\\
                          \vdots&\ddots&\vdots&\vdots\\
                         t_{i-1\ts 1}(u)&\dots&t_{i-1\ts i-1}(u)&t_{i-1\ts i}(u)\\
                         t_{j\ts 1}(u)&\dots&t_{j\ts i-1}(u)&\boxed{t_{j\tss i}(u)}\\
           \end{vmatrix}\ts h_i(u)^{-1}
\eeq
for $1\leqslant i<j\leqslant 1'$.
By \cite[Lem.~4.1]{m:dt}, under the anti-automorphism $\tau$
of $\X(\osp_{1|2m})$ defined in
\eqref{tauanti}, for all $k$ and $i<j$ we have
\beql{taue}
\tau: h_k(u)\mapsto h_k(u)\fand
e_{ij}(u)\mapsto f_{ji}(u)(-1)^{\bi\bj+\bj},\quad f_{ji}(u)\mapsto e_{ij}(u)(-1)^{\bi\bj+\bi}.
\eeq

Introduce the coefficients of the series defined in
\eqref{hmqua}, \eqref{eijmlqua} and \eqref{fijlmqua} by the expansions
\beql{enise}
e_{ij}(u)=\sum_{r=1}^{\infty} e_{ij}^{(r)}\tss u^{-r},\qquad
f_{ji}(u)=\sum_{r=1}^{\infty} f_{ji}^{(r)}\tss u^{-r},\qquad
h_i(u)=1+\sum_{r=1}^{\infty} h_i^{(r)}\tss u^{-r}.
\eeq
Furthermore, set
\beql{defkn}
k_{i}(u)=h_i(u)^{-1}h_{i+1}(u),\qquad
e_{i}(u)=e_{i\ts i+1}(u),
\qquad f_{i}(u)=f_{i+1\ts i}(u),
\eeq
for $i=1,\dots, m$.
We will also use the coefficients of the series defined by
\beql{efexp}
e_i(u)=\sum_{r=1}^{\infty}e_i^{(r)}u^{-r}\Fand
f_i(u)=\sum_{r=1}^{\infty}f_i^{(r)}u^{-r}.
\eeq

By \cite[Prop.~5.1]{m:dt}, the Gaussian generators $h_i(u)$ satisfy the relations
\beql{ilm}
h_i(u)\ts h_{i'}\big(u+m-i+1/2\big)
=h_{i+1}(u)\ts h_{(i+1)'}\big(u+m-i+1/2\big)
\eeq
for $i=1,\dots,m$. Together with the relation
\beql{cuhh}
c(u)=h_1(u)\tss h_{1'}(u+m+1/2)
\eeq
for the central series $c(u)$ defined in \eqref{ttra},
they imply that the coefficients of
all series $h_i(u)$ with $i=1,2,\dots,1'$ pairwise commute in
$\X(\osp_{1|2m})$; see \cite[Cor.~5.2]{m:dt}.
We will also recall a formula for $c(u)$ in terms
of the Gaussian generators $h_i(u)$ with $i=1,\dots,m+1$; see \cite[Thm~5.3]{m:dt}.
We have
\beql{cu}
c(u)=\prod_{i=1}^m \ts\frac{h_i(u+i-1)}{h_i(u+i)}\cdot h_{m+1}(u+m+1/2)\ts h_{m+1}(u+m).
\eeq

\section{Drinfeld-type presentations of the Yangians for $\osp_{1|2}$}
\label{sec:dpot}

We will now suppose that $m=1$ and give Drinfeld-type presentations
of the algebras $\X(\osp_{1|2})$ and $\Y(\osp_{1|2})$. Our approach is similar to
\cite{aacfr:sy}, but we use a
different set of generators by
adjoining the coefficients of the series $e_{11'}(u)$ and $f_{1'1}(u)$.
This allows us to avoid Serre-type relations used therein.
We use notation \eqref{defkn} and set $e(u)=e_1(u)$, $f(u)=f_1(u)$ and $k(u)=k_1(u)$.

\bth\label{thm:odp}
The extended Yangian $\X(\osp_{1|2})$ is generated by
the coefficients of the series
$h_1(u), h_2(u), e(u), f(u), e_{11'}(u)$ and $f_{1'1}(u)$,
subject only to the following relations.
We have
\begin{align}
\label{ohihj}
\big[h_i(u),h_j(v)\big]&=0\qquad\text{for all}\quad i,j\in\{1,2\}, \\[0.4em]
\label{oeifj}
\big[e(u),f(v)\big]&=\frac{k(u)-k(v)}{u-v}.
\end{align}
Furthermore,
\begin{align}
\label{ohiej}
\big[h_1(u),e(v)\big]&=
h_1(u)\ts\frac{e(u)-e(v)}{u-v},\\[0.4em]
\label{ohifj}
\big[h_1(u),f(v)\big]&=-
\frac{f(u)-f(v)}{u-v}\ts h_1(u)
\end{align}
and
\begin{align}
\label{ohtej}
\big[h_2(u),e(v)\big]&
=h_2(u)\,\Big(\frac{e(u) -e(v) }{u-v}
-\frac{e(u-1/2)-e(v)}{u-v-1/2}\Big),\\[0.4em]
\label{ohtfj}
\big[h_2(u),f(v)\big]&=\Big({-}\frac{f(u) -f(v) }{u-v}
+\frac{f(u-1/2)-f(v)}{u-v-1/2}\Big)\,h_2(u).
\end{align}
We also have
\begin{align}
\non
\big[e(u),e(v)\big]&=\frac{e(u)^2+e_{11'}(u)-e(v)^2-e_{11'}(v)}{u-v}\\
{}&+\frac{e(u)\tss e(v)-e(v)\tss e(u)}{2\tss(u-v)}
-\frac{\big(e(u)-e(v)\big)^2}{2\tss(u-v)^2}
\label{oeiei}
\end{align}
and
\begin{align}
\non
\big[f(u),f(v)\big]&=\frac{f(u)^2-f_{1'1}(u)-f(v)^2+f_{1'1}(v)}{u-v}\\
{}&-\frac{f(u)\tss f(v)-f(v)\tss f(u)}{2\tss(u-v)}
-\frac{\big(f(u)-f(v)\big)^2}{2\tss(u-v)^2}.
\label{ofifi}
\end{align}
Finally,
\begin{align}
\non
\big[e(u),e_{11'}(v)\big]&=-\frac{\big(e(u)-e(v)\big)
\big(e_{11'}(u)-e_{11'}(v)\big)}{u-v}\\[0.4em]
{}&-\frac{e(u+1/2)-e(v)}{u-v+1/2}\ts e(u)^2
-\frac{e_{11'}(u+1/2)-e_{11'}(v)}{u-v+1/2}\ts e(u)
\label{oeieoo}
\end{align}
and
\begin{align}
\non
\big[f(u),f_{1'1}(v)\big]&=\frac{\big(f_{1'1}(u)-f_{1'1}(v)\big)
\big(f(u)-f(v)\big)}{u-v}\\[0.4em]
{}&-f(u)^2\ts\ts\frac{f(u+1/2)-f(v)}{u-v+1/2}
+f(u)\ts\frac{f_{1'1}(u+1/2)-f_{1'1}(v)}{u-v+1/2}.
\label{ofifoo}
\end{align}
\eth

\bpf
As the first step, we will verify that all the above relations hold in the extended Yangian.
Relations \eqref{ohihj} and \eqref{oeifj} were pointed out in \cite{aacfr:sy} and
\cite[Sec.~3]{m:ry} along with the identities
\beql{ef}
e_{21'}(u)=-e(u-1/2)\Fand f_{1'2}(u)=f(u-1/2).
\eeq
It is sufficient to verify \eqref{ohiej}, \eqref{ohtej}, \eqref{oeiei} and \eqref{oeieoo},
because the remaining relations will follow by the application
of the anti-automorphism $\tau$ using \eqref{taue}. By \eqref{defrel} we have
\ben
\big[\tss t_{11}(u),t_{12}(v)\big]=-\frac{1}{u-v}
\big(t_{11}(u)\ts t_{12}(v)-t_{11}(v)\ts t_{12}(u)\big).
\een
Since $h_1(u)=t_{11}(u)$ and $e(v)=t_{11}(v)^{-1}t_{12}(v)$, by multiplying both sides
by $t_{11}(v)^{-1}$ from the left we get \eqref{ohiej}.
Furthermore, by \eqref{ilm} and \eqref{cuhh} we have
\beql{hhhc}
h_1(u)\ts h_{1'}(u+1/2)=
h_{2}(u)\ts h_{2}(u+1/2)\Fand h_1(u)\ts h_{1'}(u+3/2)=c(u).
\eeq
There exists a unique power series $z(u)$ in $u^{-1}$ with coefficients
in the center of $\X(\osp_{1|2})$ and with the constant term $1$, satisfying the relation
$z(u)\tss z(u+1/2)=c(u-1)$. This implies that $h_2(u)$ can be expressed by
\beql{htz}
h_2(u)=z(u)\tss h_1(u-1/2)\tss h_1(u-1)^{-1}.
\eeq
We will use this relation to derive \eqref{ohtej} from \eqref{ohiej}.
By rearranging the latter we get
\beql{evho}
e(v)\tss h_1(u)=h_1(u)\Big(\frac{u-v+1}{u-v}\ts e(v)-\frac{1}{u-v}\ts e(u)\Big).
\eeq
In particular, setting $v=u+1$ yields
\beql{euho}
e(u+1)\tss  h_1(u)=h_1(u)\tss e(u).
\eeq
Therefore, we have
\ben
e(v)\tss h_1(u)=\frac{u-v+1}{u-v}\ts h_1(u)\tss e(v)-\frac{1}{u-v}\ts e(u+1)\tss h_1(u)
\een
which implies
\beql{evhoinv}
e(v)\tss h_1(u)^{-1}=h_1(u)^{-1}\Big(\frac{u-v}{u-v+1}\ts e(v)+\frac{1}{u-v+1}\ts e(u+1)\Big).
\eeq
Since the series $z(u)$ is central, by using \eqref{htz} together with \eqref{evho}
and \eqref{evhoinv}, we derive the relation
\ben
e(v)\tss h_2(u)=h_2(u)\Big(\frac{(u-v+1/2)(u-v-1)}{(u-v-1/2)(u-v)}\ts e(v)+\frac{1}{u-v-1/2}\ts e(u-1/2)
-\frac{1}{u-v}\ts e(u)\Big)
\een
which is equivalent to \eqref{ohtej}.

Now consider two particular cases of \eqref{defrel},
\begin{align}
\big[\tss t_{11}(u),t_{11'}(v)\big]&=-\frac{1}{u-v}
\big(t_{11}(u)\ts t_{11'}(v)-t_{11}(v)\ts t_{11'}(u)\big)
\non\\
{}&-\frac{1}{u-v+3/2}
\big(t_{11'}(v)\ts t_{11}(u)+t_{12}(v)\ts t_{12}(u)-t_{11}(v)\ts t_{11'}(u)\big)
\label{tootoo}
\end{align}
and
\begin{align}
\big[\tss t_{12}(u),t_{12}(v)\big]&=-\frac{1}{u-v}
\big(t_{12}(u)\ts t_{12}(v)-t_{12}(v)\ts t_{12}(u)\big)
\non\\
{}&-\frac{1}{u-v+3/2}
\big(t_{11'}(v)\ts t_{11}(u)+t_{12}(v)\ts t_{12}(u)-t_{11}(v)\ts t_{11'}(u)\big).
\non
\end{align}
By expanding the super-commutators and eliminating the product $t_{11'}(v)\ts t_{11}(u)$
in the second formula using the first,
we come to the relation
\ben
-t_{11}(u)\ts t_{11'}(v)+t_{11}(v)\ts t_{11'}(u)
=(u-v+1/2)\ts t_{12}(u)\ts t_{12}(v)+(u-v-1/2)\ts t_{12}(v)\ts t_{12}(u).
\een
The right hand side equals
\ben
(u-v+1/2)\ts h_1(u)\tss e(u)\tss h_1(v)\tss e(v)+
(u-v-1/2)\ts h_1(v)\tss e(v)\tss h_1(u)\tss e(u).
\een
Transform it by applying \eqref{evho} to the products $e(u)\tss h_1(v)$ and
$e(v)\tss h_1(u)$. By taking into account $t_{11'}(u)=h_1(u)\tss e_{11'}(u)$
and multiplying from the left by the inverse of $h_1(u)h_1(v)$, we then obtain
\begin{align}
e_{11'}(u)-e_{11'}(v)&=\frac{(u-v+1/2)(u-v-1)}{u-v}\ts e(u)\ts e(v)
+\frac{(u-v-1/2)(u-v+1)}{u-v}\ts e(v)\ts e(u)
\non\\
{}&-\frac{u-v-1/2}{u-v}\ts e(u)^2+\frac{u-v+1/2}{u-v}\ts e(v)^2,
\label{eooeoo}
\end{align}
which is equivalent to \eqref{oeiei}.

Finally, to prove \eqref{oeieoo}, begin with the following particular case of \eqref{defrel},
\beql{ttoo}
\big[\tss t_{12}(u),t_{11'}(v)\big]=-\frac{1}{u-v}
\big(t_{12}(u)\ts t_{11'}(v)-t_{12}(v)\ts t_{11'}(u)\big).
\eeq
Note its consequence $t_{12}(u+1)\ts t_{11'}(u)=t_{11'}(u+1)\ts t_{12}(u)$
which implies
\beql{hehe}
h_1(u+1)\tss e(u+1)\tss h_1(u)\tss e_{11'}(u)=h_1(u+1)\tss e_{11'}(u+1)\tss h_1(u)\tss e(u).
\eeq
Write \eqref{ttoo}
in terms of the Gaussian generators and multiply both sides by $h_1(u)^{-1}h_1(v)^{-1}$
from the left to get
\begin{align}
\frac{u-v+1}{u-v}\ts h_1(v)^{-1}\tss e(u)\tss h_1(v)\tss e_{11'}(v)
{}&-h_1(u)^{-1}\tss e_{11'}(v)\tss h_1(u)\tss e(u)\non\\
{}&=\frac{1}{u-v}\ts h_1(u)^{-1}\tss e(v)\tss h_1(u)\tss e_{11'}(u).
\label{hieie}
\end{align}
Similarly, by multiplying both sides of \eqref{tootoo} by $h_1(u)^{-1}h_1(v)^{-1}$
from the left, we obtain
\ben
\bal
e_{11'}(v)&-h_1(u)^{-1}e_{11'}(v)h_1(u)=-\frac{1}{u-v}\big(e_{11'}(v)-e_{11'}(u)\big)\\[0.4em]
{}&-\frac{1}{u-v+3/2}\big(h_1(u)^{-1}e_{11'}(v)h_1(u)+h_1(u)^{-1}e(v)h_1(u)\tss e(u)-e_{11'}(u)\big).
\eal
\een
Replacing the product $e(v)h_1(u)$ by \eqref{evho} and rearranging, we come to
\begin{align}
h_1(u)^{-1}e_{11'}(v)h_1(u)&=\frac{(u-v+3/2)(u-v+1)}{(u-v)(u-v+1/2)}\ts e_{11'}(v)-
\frac{2u-2v+3/2}{(u-v)(u-v+1/2)}\ts e_{11'}(u)
\non\\[0.5em]
{}&+\frac{1}{u-v+1/2}\ts\Big(\frac{u-v+1}{u-v}\ts e(v)-\frac{1}{u-v}\ts e(u)\Big)\ts e(u).
\label{heh}
\end{align}
Substitute this expression into \eqref{hieie} and
apply \eqref{evho} to the products $e(u)\tss h_1(v)$ and
$e(v)\tss h_1(u)$. Multiplying both sides by $(u-v)/(u-v+1)$, we come to the relation
\begin{multline}
\non
\big[e(u),e_{11'}(v)\big]=\frac{e(u)-e(v)}{u-v}\ts e_{11'}(v)-
\Big(\frac{e(u)}{(u-v)(u-v+1)}-\frac{e(v)}{u-v}\Big)\ts e_{11'}(u)\\[0.4em]
{}+\frac{1}{u-v+1/2}\Big(e_{11'}(v)-\frac{2u-2v+3/2}{u-v+1}\ts e_{11'}(u)\Big)\ts e(u)
\\[0.4em]
{}+\frac{1}{u-v+1/2}\Big(e(v)-\frac{1}{u-v+1}\ts e(u)\Big)\ts e(u)^2.
\end{multline}
On the other hand, setting $v=u+1$ into \eqref{heh}, we get
\ben
e_{11'}(u+1)\tss h_1(u)=h_1(u)\tss\big(e_{11'}(u)-2\tss e(u)^2\big).
\een
Together with \eqref{euho} and \eqref{hehe} this yields
\beql{eueoo}
\big[e(u),e_{11'}(u)\big]=-2\tss e(u)^3.
\eeq
By using
this identity we can simplify the above formula for the super-commutator to
\begin{multline}
\non
\big[e(u),e_{11'}(v)\big]=\frac{e(u)-e(v)}{u-v}\ts \big(e_{11'}(v)-e_{11'}(u)\big)\\[0.4em]
{}+\frac{1}{u-v+1/2}\Big(\big(e_{11'}(v)-e_{11'}(u)\big)\ts e(u)+e(v)\tss e(u)^2-2\tss e(u)^3\Big).
\end{multline}
Set $v=u+1/2$ in \eqref{eooeoo} to get another identity
\beql{eoomeoo}
e_{11'}(u+1/2)-e_{11'}(u)+e(u+1/2)\tss e(u)-2\tss e(u)^2=0.
\eeq
Its use brings the above relation for $[e(u),e_{11'}(v)]$
to the required form \eqref{oeieoo}.

Since all relations in the formulation of the theorem hold in the extended Yangian, we have
a homomorphism
\beql{homexy}
\wh \X(\osp_{1|2})\to\X(\osp_{1|2}),
\eeq
where $\wh \X(\osp_{1|2})$ denotes the algebra
whose (abstract) generators are the
coefficients of series
$h_1(u), h_2(u), e(u), f(u), e_{11'}(u)$ and $f_{1'1}(u)$
given by the same expansions as in \eqref{enise} and \eqref{efexp},
with the relations as in the statement of the theorem
(omitting the subscripts of $e_{12}$ and $f_{21}$).
The homomorphism \eqref{homexy} takes the generators to the elements
of $\X(\osp_{1|2})$ with the same name. We will show that
this homomorphism is surjective and injective. The surjectivity is clear
from the Gauss decomposition \eqref{gd}, formulas \eqref{ef}
and the first relation in \eqref{hhhc}.

Now we prove the injectivity of the homomorphism \eqref{homexy}.
The same application of the Poincar\'e--Birkhoff--Witt theorem for the algebra
$\X(\osp_{1|2})$ as in \cite[Sec.~6]{m:dt} shows that
the set of monomials in
the generators
$h_{1}^{(r)}, h_{2}^{(r)},e^{(r)},f^{(r)},e_{11'}^{(r)}$ and $f_{1'1}^{(r)}$
with $r\geqslant 1$
taken in some fixed order, with the powers of odd generators
not exceeding $1$,
is linearly independent in the extended Yangian
$\X(\osp_{1|2})$. Therefore, to complete the proof of the theorem, it is sufficient
to verify that the monomials in the generators
$h_{1}^{(r)}, h_{2}^{(r)},e^{(r)},f^{(r)},e_{11'}^{(r)}$ and $f_{1'1}^{(r)}$
with $r\geqslant 1$ of the algebra $\wh \X(\osp_{1|2})$, taken in a certain
fixed order, span the algebra.

Define the ascending filtration
on the algebra $\wh \X(\osp_{1|2})$ by setting the degree of each generator
with the superscript $r$ to be equal to $r-1$. We will use the bar symbol to
denote the image of each generator in the
$(r-1)$-th component of the graded algebra $\gr\wh \X(\osp_{1|2})$.
The defining relations of $\wh \X(\osp_{1|2})$ imply the corresponding relations
for these images in the graded algebra, which are easily derived with the use of the expansion
formula
\beql{expafo}
\frac{g(u)-g(v)}{u-v}=-\sum_{r,s\geqslant 1}\tss g^{(r+s-1)}\tss u^{-r}v^{-s}\qquad\text{for}\quad
g(u)=\sum_{k=1}^{\infty}\tss g^{(k)}\tss u^{-k}.
\eeq
Namely, relations \eqref{ohihj} -- \eqref{ohtfj} imply
\ben
\big[\hba^{(r)}_i,\hba^{(s)}_j\big]=0,\quad
\big[\eb^{(r)},\fb^{(s)}\big]=\hba_1^{(r+s-1)}-\hba_2^{(r+s-1)}
\een
and
\ben
\big[\hba^{(r)}_1,\eb^{(s)}\big]=-\eb^{(r+s-1)},\quad
\big[\hba^{(r)}_1,\fb^{(s)}\big]=\fb^{(r+s-1)},
\quad
\big[\hba^{(r)}_2,\eb^{(s)}\big]=
\big[\hba^{(r)}_2,\fb^{(s)}\big]=0,
\een
while relations \eqref{oeiei} -- \eqref{ofifoo} give
\ben
\big[\eb^{(r)},\eb^{(s)}\big]=-\eb_{11'}^{(r+s-1)},\quad
\big[\fb^{(r)},\fb^{(s)}\big]=\fb_{1'1}^{(r+s-1)},
\quad
\big[\eb^{(r)},\eb_{11'}^{(s)}\big]=
\big[\fb^{(r)},\fb_{1'1}^{(s)}\big]=0.
\een
This determines all
super-commutator relations between the generators of $\gr\wh \X(\osp_{1|2})$.
In particular, we have
\ben
\big[\eb^{(r)},\fb_{1'1}^{(s)}\big]=2\tss \fb^{(r+s-1)},\quad
\big[\eb_{11'}^{(r)},\fb^{(s)}\big]=-2\tss \eb^{(r+s-1)},\quad
\big[\eb_{11'}^{(r)},\fb_{1'1}^{(s)}\big]=4\tss (\hba_2^{(r+s-1)}-\hba_1^{(r+s-1)}).
\een
The spanning property of the ordered monomials now follows from the observation that
the super-commutator relations coincide with those in the polynomial current Lie
superalgebra $\agot[u]$, where $\agot$ is the centrally extended Lie superalgebra
$\osp_{1|2}$. This
completes the proof of the theorem.
\epf

The following is a version of the Poincar\'e--Birkhoff--Witt
theorem for the orthosymplectic Yangian which was established in the proof
of Theorem~\ref{thm:odp}.

\bco\label{cor:opbwdp}
The set of monomials in the elements
$h_{1}^{(r)}, h_{2}^{(r)},e^{(r)},f^{(r)},e_{11'}^{(r)}$ and $f_{1'1}^{(r)}$,
where $r=1,2,\dots$,
taken in some fixed order, with the powers of odd generators
not exceeding $1$,
forms a basis of $\X(\osp_{1|2})$.
\qed
\eco

By the definition of the Gaussian generators, the coefficients of all series
$k(u)$, $e(u)$, $f(u)$, $e_{11'}(u)$ and $f_{1'1}(u)$ are stable under the action
of all automorphisms \eqref{muf} and so they
belong to the subalgebra $\Y(\osp_{1|2})$
of $\X(\osp_{1|2})$. We now derive
a Drinfeld-type presentation
of the Yangian $\Y(\osp_{1|2})$.

\bco\label{cor:odpy}
The Yangian $\Y(\osp_{1|2})$ is generated by
the coefficients of the series
$k(u)$, $e(u)$, $f(u)$, $e_{11'}(u)$ and $f_{1'1}(u)$,
subject only to relations \eqref{oeifj},
\eqref{oeiei} -- \eqref{ofifoo} together with
\beql{kukv}
\big[k(u),k(v)\big]=0,
\eeq
\beql{kuev}
\big[k(u),e(v)\big]
=k(u)\ts\Big({-}\frac{e(u-1/2) -e(v) }{3\tss(u-v-1/2)}
-\frac{2\tss\big(e(u+1)-e(v)\big)}{3\tss(u-v+1)}\Big)
\eeq
and
\beql{kufv}
\big[k(u),f(v)\big]=\Big(\frac{f(u-1/2) -f(v) }{3\tss(u-v-1/2)}
+\frac{2\tss\big(f(u+1)-f(v)\big)}{3\tss(u-v+1)}\Big)\ts k(u).
\eeq
\eco

\bpf
Relation \eqref{kukv} follows from \eqref{ohihj}, so we only need to verify
\eqref{kuev}, because \eqref{kufv} will then follow by the application of
the anti-automorphism $\tau$ via \eqref{taue}. Since $h_2(u)=k(u)\tss h_1(u)$,
we can write
\ben
\big[h_2(u),e(v)\big]=k(u)\ts \big[h_1(u),e(v)\big]+\big[k(u),e(v)\big]\ts h_1(u)
\een
and so
\ben
\big[k(u),e(v)\big]=\big[h_2(u),e(v)\big]\ts h_1(u)^{-1}
-h_2(u)\ts h_1(u)^{-1}\big[h_1(u),e(v)\big]\ts h_1(u)^{-1}.
\een
Now apply \eqref{ohiej} and \eqref{ohtej} to the super-commutators on the right hand side
to get
\ben
\big[k(u),e(v)\big]=
-k(u)\tss h_1(u)\ts \frac{e(u-1/2)-e(v)}{u-v-1/2}\ts h_1(u)^{-1}.
\een
The derivation of \eqref{kuev} is completed by the application of the following
consequence of \eqref{evhoinv},
\ben
h_1(u)\ts e(v)\tss h_1(u)^{-1}=\frac{u-v}{u-v+1}\ts e(v)+\frac{1}{u-v+1}\ts e(u+1).
\een

It is clear from the decomposition \eqref{tensordecom} (with $m=1$) that the coefficients
of the series $k(u)$, $e(u)$, $f(u)$, $e_{11'}(u)$ and $f_{1'1}(u)$ generate
the subalgebra $\Y(\osp_{1|2})$
of the extended Yangian $\X(\osp_{1|2})$; cf. \cite[Prop.~6.1]{jlm:ib}. Therefore,
we have an epimorphism from the (abstract) algebra $\wh\Y(\osp_{1|2})$
defined by the generators and relations as in the statement of the corollary,
to the Yangian $\Y(\osp_{1|2})$, which
takes the generators to the elements
of $\Y(\osp_{1|2})$ denoted by the same symbols. Given any series
$\vp(u)\in 1+u^{-1}\CC[[u^{-1}]]$, consider the automorphism
of the algebra $\wh\X(\osp_{1|2})$ introduced in the proof of Theorem~\ref{thm:odp},
defined by
\ben
h_i(u)\mapsto \vp(u)\tss h_i(u)\qquad\text{for}\quad i=1,2,
\een
and which leaves all the remaining generators fixed; cf.~\eqref{muf}.  Then
$\wh \Y(\osp_{1|2})$ coincides with the subalgebra of $\wh\X(\osp_{1|2})$
which consists of the elements stable under all these automorphisms.
Therefore, the
epimorphism $\wh \Y(\osp_{1|2})\to \Y(\osp_{1|2})$ can be regarded as the restriction
of the isomorphism $\wh \X(\osp_{1|2})\to \X(\osp_{1|2})$, and hence is injective.
\epf

\bco\label{cor:opbwdpy}
The set of monomials in the elements
$k^{(r)},e^{(r)},f^{(r)},e_{11'}^{(r)}$ and $f_{1'1}^{(r)}$,
where $r=1,2,\dots$,
taken in some fixed order, with the powers of odd generators
not exceeding $1$,
forms a basis of $\Y(\osp_{1|2})$.
\qed
\eco

By taking the coefficients of $v^0$ on both sides of \eqref{eooeoo},
and applying \eqref{taue}, we get
\beql{eeff}
e_{11'}(u)=-e(u)^2-[e^{(1)},e(u)]\Fand f_{1'1}(u)=f(u)^2+[f^{(1)},f(u)].
\eeq
Therefore, the coefficients of the series $e_{11'}(u)$ and $f_{1'1}(u)$
can be eliminated from the Yangian defining relations.
In other words, we may regard the Yangian $\Y(\osp_{1|2})$ as the algebra
with generators $k^{(r)}, e^{(r)}$ and $f^{(r)}$ subject to the relations
of Corollary~\ref{cor:odpy}, where all occurrences of $e_{11'}(u)$ and $f_{1'1}(u)$
are replaced by \eqref{eeff}. This was the viewpoint taken in \cite{aacfr:sy},
where a different presentation of $\Y(\osp_{1|2})$ was given with the use of
certain Serre-type
relations.

To make a more explicit connection with the presentation of the Yangian $\Y(\osp_{1|2})$
given in
\cite[Theorem~3.1]{aacfr:sy}, we will use the automorphism $\si$ defined in \eqref{sigma}.
Observe that the subalgebra $\Y(\osp_{1|2})$ of
$\X(\osp_{1|2})$ is stable under $\si$.
We will keep the same notation for the restriction of $\si$ to $\Y(\osp_{1|2})$.

\ble\label{lem:sigauss}
The images of the generators of the algebra $\Y(\osp_{1|2})$ under the automorphism $\si$
are given by
\ben
\si:k(u)\mapsto k(-u),\qquad e(u)\mapsto f(-u),\qquad f(u)\mapsto -e(-u).
\een
\ele

\bpf
Since $e(u)=t_{11}(u)^{-1}\tss t_{12}(u)$, we find
\ben
\si:e(u)\mapsto t_{11}(-u-1)^{-1}\tss t_{21}(-u-1)=t_{21}(-u)\ts t_{11}(-u)^{-1}=f(-u),
\een
where we used the relation $t_{11}(u)\ts t_{21}(u+1)=t_{21}(u)\ts t_{11}(u+1)$ implied by \eqref{defrel}.
Similarly,
\ben
\si:f(u)\mapsto -t_{12}(-u-1)\ts t_{11}(-u-1)^{-1}=-t_{11}(-u)^{-1}\ts t_{12}(-u)=-e(-u).
\een
To calculate the image of $k(u)$, first find the image of the series $c(u)$ defined in \eqref{ttra}.
By taking the $(1,1)$-entry of the first matrix product in \eqref{ttra}, we get
\ben
c(u)=t_{11}(u+3/2)\ts t_{1'1'}(u)-t_{12}(u+3/2)\ts t_{1'2}(u)-t_{11'}(u+3/2)\ts t_{1'1}(u).
\een
Hence, the image of $c(u)$ under the map $\si$ equals
\ben
t_{11}(-u-5/2)\ts t_{1'1'}(-u-1)-t_{21}(-u-5/2)\ts t_{21'}(-u-1)-t_{1'1}(-u-5/2)\ts t_{11'}(-u-1).
\een
Therefore, $\si:c(u)\mapsto c(-u-5/2)$ which follows by taking the $(1',1')$-entry of
the second matrix product in \eqref{ttra}. We can now find the image of the series $h_2(u)$ by using \eqref{htz}.
Since $c(u)=z(u+1)\tss z(u+3/2)$ and the series $z(u)$ is uniquely determined by this relation,
we derive that $\si:z(u)\mapsto z(-u)$ and so
\ben
\si:h_2(u)\mapsto z(-u)\tss h_1(-u-1/2)\tss h_1(-u)^{-1}=h_2(-u)\tss h_1(-u-1)\tss h_1(-u)^{-1}.
\een
This implies that $\si:k(u)\mapsto k(-u)$.
\epf

The following corollary essentially reproduces \cite[Theorem~3.1]{aacfr:sy} (in our notation).

\bco\label{cor:serre}
The Yangian $\Y(\osp_{1|2})$ is generated by
the coefficients of the series
$k(u)$, $e(u)$ and $f(u)$
subject only to relations \eqref{oeifj},
\eqref{kukv} and \eqref{kuev} together with
\beql{akufv}
\big[k(u),f(v)\big]=k(u)\Big(\frac{f(u+1/2) -f(v) }{3\tss(u-v+1/2)}
+\frac{2\tss\big(f(u-1)-f(v)\big)}{3\tss(u-v-1)}\Big),
\eeq
\beql{woeiei}
\big[e(u),e(v)\big]=-\frac{[e^{(1)},e(u)-e(v)]}{u-v}
+\frac{e(u)\tss e(v)-e(v)\tss e(u)}{2\tss(u-v)}
-\frac{\big(e(u)-e(v)\big)^2}{2\tss(u-v)^2},
\eeq
\beql{wofifi}
\big[f(u),f(v)\big]=-\frac{[f^{(1)},f(u)-f(v)]}{u-v}
-\frac{f(u)\tss f(v)-f(v)\tss f(u)}{2\tss(u-v)}
-\frac{\big(f(u)-f(v)\big)^2}{2\tss(u-v)^2},
\eeq
and the Serre-type relations
\begin{align}\label{serreacfr}
e(u)^3&=e(u)\ts [e(u),e^{(1)}]+[e^{(1)\tss 2},e(u)],\\[0.4em]
f(u)^3&=-f(u)\ts[f(u),f^{(1)}]+[f^{(1)\tss 2},f(u)].
\label{aserreacfrf}
\end{align}
\eco

\bpf
First we verify that all relations hold in the algebra $\Y(\osp_{1|2})$.
Relation \eqref{akufv} follows by the application of
the automorphism $\si$ to both sides of \eqref{kuev}.
It is clear that \eqref{woeiei} and \eqref{wofifi} are immediate from
\eqref{oeiei} and \eqref{ofifi} due to \eqref{eeff}. Furthermore,
by taking the coefficient of $v^{-1}$ in \eqref{oeieoo} we get
\ben
[e(u),e_{11'}^{(1)}]=
e(u)\ts e_{11'}(u)+e(u+1/2)\ts e(u)^2
+e_{11'}(u+1/2)\ts e(u).
\een
Now use \eqref{eoomeoo} to write this in the form
\ben
[e(u),e_{11'}^{(1)}]=
e(u)\ts e_{11'}(u)+e_{11'}(u)\ts e(u)
+2\ts e(u)^3
\een
which together with \eqref{eueoo} (a consequence of \eqref{oeiei} and \eqref{oeieoo})
yield
\ben
[e(u),e_{11'}^{(1)}]=
2\ts e(u)\ts e_{11'}(u)+4\ts e(u)^3.
\een
It remains to replace $e_{11'}(u)$ by \eqref{eeff}
and note that $e_{11'}^{(1)}=-2\tss e^{(1)\tss 2}$
to arrive at \eqref{serreacfr}; relation \eqref{aserreacfrf} then follows by the application of $\si$.

We thus have
an epimorphism
\beql{yhomexy}
\wh \Y(\osp_{1|2})\to\Y(\osp_{1|2}),
\eeq
where $\wh \Y(\osp_{1|2})$ denotes the algebra
whose (abstract) generators are the
coefficients of series
$k(u), e(u)$ and $f(u)$
with the relations as in the statement of the corollary.
The epimorphism \eqref{yhomexy} takes the generators to the elements
of $\Y(\osp_{1|2})$ with the same name. We only need to show that it is injective.
Introduce the series $e_{11'}(u)$ and $f_{1'1}(u)$ with coefficients
in the algebra $\wh \Y(\osp_{1|2})$ by formulas \eqref{eeff} and proceed
as in the proof of Theorem~\ref{thm:odp}. It is sufficient to show that
the monomials in the generators
$k^{(r)},e^{(r)},f^{(r)},e_{11'}^{(r)}$ and $f_{1'1}^{(r)}$
with $r\geqslant 1$ of the algebra $\wh \Y(\osp_{1|2})$, taken in a certain
fixed order, span the algebra.

Define the ascending filtration
on the algebra $\wh \Y(\osp_{1|2})$ by setting the degree of each generator
with the superscript $r$ to be equal to $r-1$ and use the bar symbol to
denote the image of each generator in the
$(r-1)$-th component of the associated graded algebra $\gr\wh \Y(\osp_{1|2})$.

Relations \eqref{eeff} and \eqref{woeiei} imply $[\eb^{(r)},\eb^{(s)}]=-\eb_{11'}^{(r+s-1)}$,
and we derive from \eqref{kuev}
that $[\bk^{(2)},\eb^{(s)}]=\eb^{(s+1)}$. Hence,
$[\bk^{(2)},\eb_{11'}^{(s)}]=2\tss \eb_{11'}^{(s+1)}$.
Now write the Serre-type relation
\eqref{serreacfr} in a different form. Note that
\ben
[e^{(1)\tss 2},e(u)]=\big[e^{(1)},[e^{(1)},e(u)]\big]-\big[[e^{(1)},e(u)],e^{(1)}\big].
\een
Therefore, by using the relation
\ben
[e^{(1)},e(u)]=-e_{11'}(u)-e(u)^2
\een
implied by \eqref{eeff}, we can write \eqref{serreacfr} as
\ben
\big[e^{(1)},e_{11'}(u)\big]=-2\tss e(u)^3-e(u)\tss e_{11'}(u)-[e^{(1)},e(u)^2],
\een
which yields $\big[\eb^{(1)},\eb_{11'}^{(s)}\big]=0$. Therefore, by taking repeated
commutators with $\bk^{(2)}$ we obtain $[\eb^{(r)},\eb_{11'}^{(s)}]=0$.
Together with their counterparts for the elements $\fb^{(r)}$, these relations
coincide with those in the polynomial current Lie
superalgebra $\osp_{1|2}[u]$ thus implying
the desired spanning property of the ordered monomials.
\epf

The generators of $\Y(\osp_{1|2})$ used in
Corollary~\ref{cor:serre} and in \cite[Theorem~3.1]{aacfr:sy} are related as follows.
The series $e(u)$ is the same,
$k(u)$ corresponds to $h(u)$ in \cite{aacfr:sy} and our $f(u)$ corresponds to $-f(u)$
in \cite{aacfr:sy}. The different-looking relations (3.3) and (3.4) in \cite{aacfr:sy}
are in fact equivalent to \eqref{kuev} and \eqref{akufv}, respectively.
Indeed, to outline the calculation, write (3.3) in our notation and rearrange to get
\begin{multline}
\big[k(u),e(v)\big]=\frac{u-v}{(u-v-1/2)(u-v+1)}\ts k(u)\ts e(v)\\[0.4em]
{}-\frac{1}{2\tss (u-v-1/2)(u-v+1)}\ts\big[k(u),e(u)\big]
-\frac{u-v}{(u-v-1/2)(u-v+1)}\ts \big[k(u),e^{(1)}\big].
\non
\end{multline}
Take the residue at $u-v=1/2$ to derive the relation
\ben
\big[k(u),e(u)\big]=k(u)\ts e(u-1/2)-\big[k(u),e^{(1)}\big]
\een
and use it to replace
the commutator $[k(u),e(u)]$ in the previous formula. Now take the residue at $u-v=-1$
in the resulting expression to get
\ben
\big[k(u),e^{(1)}\big]=2\tss k(u)\tss e(u+1)-k(u)\tss e(u-1/2).
\een
The use of this relation brings \cite[(3.3)]{aacfr:sy} to the form \eqref{kuev}.
The equivalence of \cite[(3.4)]{aacfr:sy} and \eqref{akufv} now follows by applying
the automorphism $\si$.

It is clear from the proof of Corollary~\ref{cor:serre}
that its version, where
the counterparts of relations \eqref{kuev}, \eqref{woeiei} and \eqref{serreacfr}
involving the series $f(u)$ are derived by using the anti-automorphism $\tau$ instead of
the automorphism $\si$, is also valid. In that version, relation \eqref{akufv}
is replaced by \eqref{kufv},
relation \eqref{wofifi} is unchanged, whereas the Serre-type relation \eqref{aserreacfrf}
is replaced by
\beql{serreacfrf}
f(u)^3=[f(u),f^{(1)}]\ts f(u)-[f^{(1)\tss 2},f(u)].
\eeq

Coproduct formulas in the Hopf algebra $\Y(\osp_{1|2})$ were derived in \cite{aacfr:sy}.
They can be re-written in terms of the presentation given in Corollary~\ref{cor:odpy}
as the next proposition shows.

\bpr\label{prop:copr}
For the images of the generator series under the coproduct map we have
\begin{multline}
\De:e(u)\mapsto 1\ot e(u)+\sum_{r=0}^{\infty}(-1)^r
\Big(e(u)\ot f(u+1)+e_{11'}(u)\ot\big(f_{1'1}(u+1)-2\tss f(u+1)^2\big)\Big)^r\\
{}\times
\Big(e(u)\ot 1+e_{11'}(u)\ot\big(\frac13\ts f(u-1/2)+\frac23\ts f(u+1)\big)\Big)\big(1\ot k(u)\big)
\non
\end{multline}
and
\begin{multline}
\De:f(u)\mapsto f(u)\ot 1+\big(k(u)\ot 1\big)
\Big(1\ot f(u)-\big(\frac13\ts e(u-1/2)+\frac23\ts e(u+1)\big)\ot f_{1'1}(u)\Big)\\
{}\times\sum_{r=0}^{\infty}(-1)^r
\Big(e(u+1)\ot f(u)+\big(e_{11'}(u+1)+2\tss e(u+1)^2\big)\ot f_{1'1}(u)\Big)^r.
\non
\end{multline}
\epr

\bpf
The argument is the same as in \cite{aacfr:sy}: we
write the generator series as $e(u)=t_{11}(u)^{-1}t_{12}(u)$
and $f(u)=t_{21}(u)\ts t_{11}(u)^{-1}$ and apply
definition \eqref{Delta}. To give more details for the first formula,
write $\De(e(u))$ as
\begin{multline}
\De(t_{11}(u))^{-1}\De(t_{12}(u))
=\Big(t_{11}(u)\ot t_{11}(u)+t_{12}(u)\ot t_{21}(u)+t_{11'}(u)\ot t_{1'1}(u)\Big)^{-1}\\
{}\times\Big(t_{11}(u)\ot t_{12}(u)+t_{12}(u)\ot t_{22}(u)+t_{11'}(u)\ot t_{1'2}(u)\Big)
\non
\end{multline}
which equals
\begin{multline}
\Big(1\ot 1+e(u)\ot t_{11}(u)^{-1}t_{21}(u)+e_{11'}(u)\ot t_{11}(u)^{-1}t_{1'1}(u)\Big)^{-1}\\
{}\times\Big(1\ot e(u)+e(u)\ot t_{11}(u)^{-1}t_{22}(u)+e_{11'}(u)\ot t_{11}(u)^{-1}t_{1'2}(u)\Big).
\non
\end{multline}
As we observed in the proof of Lemma~\ref{lem:sigauss}, $t_{11}(u)^{-1}t_{21}(u)=f(u+1)$
which implies the relation
$h_1(u)f(u+1)=f(u) h_1(u)$. Moreover, relation \eqref{ohifj} implies $[h_1(u),f^{(1)}]=f(u) h_1(u)$.
Hence, \eqref{eeff} yields
\ben
h_1(u)f_{1'1}(u+1)=\big(f_{1'1}(u)+2\ts f(u)^2\big)\ts h_1(u).
\een
Therefore,
\ben
t_{11}(u)^{-1}t_{1'1}(u)=h_1(u)^{-1}f_{1'1}(u)h_1(u)
=f_{1'1}(u+1)-2\ts f(u+1)^2.
\een
Furthermore, by Gauss decomposition,
\ben
t_{11}(u)^{-1}t_{22}(u)=h_1(u)^{-1}\big(h_2(u)+f(u)h_1(u)e(u)\big)=k(u)+f(u+1)e(u).
\een
Finally, use the Gauss decomposition again to write
\ben
t_{11}(u)^{-1}t_{1'2}(u)=h_1(u)^{-1}\big(f_{1'2}(u)h_2(u)+f_{1'1}(u)h_1(u)e(u)\big)
\een
and recall that $f_{1'2}(u)=f(u-1/2)$ by \eqref{ef}.
Now re-arrange \eqref{ohifj} to bring it
to the form
\ben
h_1(u)^{-1}f(v)=\frac{1}{u-v+1}\Big((u-v)f(v)+f(u+1)\Big) h_1(u)^{-1},
\een
which implies
\ben
h_1(u)^{-1}f_{1'2}(u)=h_1(u)^{-1}f(u-1/2)
=\Big(\frac13\ts f(u-1/2)+\frac23\ts f(u+1)\Big)\ts h_1(u)^{-1}.
\een
The required formula for $\De(e(u))$ follows by expressing this image in terms of the generators used in
Corollary~\ref{cor:odpy}.
\epf

The image of the series $k(u)$ under the coproduct $\De$ can be found by using the relations
$k(u)=[e^{(1)},f(u)]=[e(u),f^{(1)}]$ implied by \eqref{oeifj}, although its explicit
expression has a rather complicated form.

\section{Yangian presentations}
\label{sec:dp}

By the embedding theorem \cite[Thm~3.1]{m:dt}, the extended Yangian
$\X(\osp_{1|2l})$ with $l<m$ can be regarded as a subalgebra of $\X(\osp_{1|2m})$.
Moreover, the embedding $\X(\osp_{1|2l})\hra \X(\osp_{1|2m})$
is consistent with the Gauss decompositions. Therefore,
the relations of Theorem~\ref{thm:odp} will hold in $\X(\osp_{1|2m})$.
We will also need the embedding for $l=2$ and first derive
some additional relations in this case.

\subsection{Relations in $\X(\osp_{1|4})$}
\label{subsec:rellt}

We will use the Gaussian generators of $\X(\osp_{1|2m})$ with $m=2$, as introduced in Sec.~\ref{sec:gd}.

\bpr\label{prop:commu}
We have the relation in $\X(\osp_{1|4})${\rm :}
\ben
\big[e_{12}^{(1)},e_{22'}(v)\big]=e_{12}(v)\tss e_{22'}(v)-e_{12'}(v)-e_{21'}(v).
\een
\epr

\bpf
By \eqref{hmqua} and \eqref{eijmlqua} we have $e_{12}^{(1)}=t_{12}^{(1)}$ along with
\ben
h_2(v)=t_{22}(v)-t_{21}(v)\ts t_{11}(v)^{-1}\tss t_{12}(v)
\een
and
\ben
e_{22'}(v)=h_2(v)^{-1}\big(t_{22'}(v)-t_{21}(v)\ts t_{11}(v)^{-1}\tss t_{12'}(v)\big).
\een
The defining relations \eqref{defrel} give
$[\tss t^{(1)}_{12},t_{11}(v)]=t_{12}(v)$ and hence
\ben
[\tss t^{(1)}_{12},t_{11}(v)^{-1}]=-t_{11}(v)^{-1}\tss t_{12}(v)\tss t_{11}(v)^{-1}.
\een
Therefore, using the commutators of $t^{(1)}_{12}$ with $t_{21}(v)$, $t_{22}(v)$ and $t_{12}(v)$
implied by \eqref{defrel}, we get
\begin{multline}
\big[\tss t^{(1)}_{12},h_2(v)\big]=-t_{12}(v)
+\big(t_{11}(v)-t_{22}(v)\big)t_{11}(v)^{-1}\tss t_{12}(v)\\[0.4em]
+t_{21}(v)\ts t_{11}(v)^{-1}\tss t_{12}(v)\ts t_{11}(v)^{-1}\tss t_{12}(v)
=-h_2(v)\ts e_{12}(v),
\non
\end{multline}
implying that $[\tss t^{(1)}_{12},h_2(v)^{-1}]=e_{12}(v)\tss h_2(v)^{-1}$.
A similar calculation yields
\ben
\big[\tss t^{(1)}_{12},t_{22'}(v)-t_{21}(v)\ts t_{11}(v)^{-1}\tss t_{12'}(v)\big]
=-h_2(v)\ts \big(e_{12'}(v)+e_{21'}(v)\big),
\een
and the required formula follows.
\epf

\bpr\label{prop:idetr}
We have the identity in $\X(\osp_{1|4})${\rm :}
\ben
e_{21'}(u)=e_{12'}(u-3/2)-e_{23}(u)\ts e_{13}(u-3/2)-e_{22'}(u)\ts e_{12}(u-3/2).
\een
\epr

\bpf
By inverting the matrices on both sides of \eqref{gd}, we get
\ben
T(u)^{-1}=E(u)^{-1}\tss H(u)^{-1}\tss F(u)^{-1}.
\een
On the other hand, relation \eqref{ttra} implies $T^{\tss t}(u+\ka)=c(u+\ka)\tss T(u)^{-1}$.
Hence, by calculating the entries of the matrix $E(u)^{-1}$
and equating the $(i,1')$ entries with $i=2,3,4$ in this matrix relation, we derive
\begin{multline}
-h_1(u+\ka)\tss e_{12'}(u+\ka)=c(u+\ka)\\
{}\times\tss \big({-}e_{21'}(u)+e_{23}(u)\tss e_{31'}(u)
+e_{22'}(u)\tss e_{2'1'}(u)-e_{23}(u)\tss e_{32'}(u)\tss e_{2'1'}(u)\big)
\tss h_{1'}(u)^{-1},
\non
\end{multline}
\ben
h_1(u+\ka)\tss e_{13}(u+\ka)=c(u+\ka)
\big({-}e_{31'}(u)+ e_{32'}(u)\tss e_{2'1'}(u)\big)
\tss h_{1'}(u)^{-1},
\een
and
\ben
h_1(u+\ka)\tss e_{12}(u+\ka)=-c(u+\ka)
\tss e_{2'1'}(u)
\tss h_{1'}(u)^{-1}.
\een
Observe that relation \eqref{ohiej} holds in the same form in $\X(\osp_{1|4})$, when
$e(u)$ is replaced with $e_{12}(u)$, $e_{13}(u)$ or $e_{12'}(u)$, thus implying
$
h_1(u)\tss e(u)=e(u+1)\tss h_1(u).
$
Furthermore,
\ben
c(u+\ka)\tss h_{1'}(u)^{-1}=h_1(u+\ka)
\een
by \eqref{cuhh}, so that replacing $\ka$ by its value $-5/2$ we derive
\ben
\bal
e_{12'}(u-3/2)&=e_{21'}(u)-e_{23}(u)\tss e_{31'}(u)
-e_{22'}(u)\tss e_{2'1'}(u)+e_{23}(u)\tss e_{32'}(u)\tss e_{2'1'}(u),\\
e_{13}(u-3/2)&=-e_{31'}(u)+ e_{32'}(u)\tss e_{2'1'}(u),\\
e_{12}(u-3/2)&=-e_{2'1'}(u),
\eal
\een
which yields the required identity.
\epf

\bco\label{cor:reid}
In the algebra $\X(\osp_{1|4})$ we have
\ben
\big[e_{12}^{(1)},e_{22'}(v)\big]=e_{12}(v)\tss e_{22'}(v)-e_{12'}(v)-
e_{12'}(v-3/2)+e_{23}(v)\ts e_{13}(v-3/2)+e_{22'}(v)\ts e_{12}(v-3/2)
\een
and
\begin{multline}
\big[e_{12}(u),e_{22'}(v)\big]=-\frac{e_{12}(u)-e_{12}(v)}{u-v}\ts e_{22'}(v)
+\frac{e_{12'}(u)-e_{12'}(v)}{u-v}+\frac{e_{12'}(u)-e_{12'}(v-3/2)}{u-v+3/2}\\[0.4em]
-e_{23}(v)\ts \frac{e_{13}(u)-e_{13}(v-3/2)}{u-v+3/2}-
e_{22'}(v)\ts \frac{e_{12}(u)-e_{12}(v-3/2)}{u-v+3/2}.
\non
\end{multline}
\eco

\bpf
The first relation is immediate from Propositions~\ref{prop:commu} and \ref{prop:idetr},
while the second follows from the first by commuting both sides with $h_1(u)$.
Here we rely on \cite[Cor.~3.3]{m:dt} implying that $h_1(u)$ commutes with each
of the series $e_{22'}(v)$ and $e_{23}(v)$, and use the commutation relation
$[h_1(u),e_{12}^{(1)}]=-h_1(u)e_{12}(u)$ which follows from \eqref{defrel}.
\epf

We point out a consequence of the second relation to be used below. By taking the coefficients
of both sides at $v^{-1}$, we get
\beql{paret}
\big[e^{(1)}_{22'},e_{12}(u)\big]=2\tss e_{12'}(u).
\eeq

\subsection{Presentations of the Yangians for $\osp_{1|2m}$}
\label{subsec:py}

Suppose that $\ve_1,\dots,\ve_{m}$ is an orthogonal basis of a vector space
with the bilinear form such that $(\ve_i,\ve_i)=-1$ for $i=1,\dots,m$.
We will take the family of vectors
\ben
\al_{i\tss j}=\ve_i-\ve_j,\qquad \al_{i\tss j'}=\ve_i+\ve_j
\qquad\text{for}\quad 1\leqslant i<j\leqslant m,
\een
and
\ben
\al_{i\ts m+1}=\ve_i,\qquad \al_{i\tss i'}=2\tss\ve_i
\qquad\text{for}\quad 1\leqslant i\leqslant m,
\een
as a system of positive roots for $\osp_{1|2m}$.
The simple roots are $\al_1,\dots,\al_{m}$ with
$
\al_i=\al_{i\ts i+1}
$
for $i=1,\dots,m$.
The associated Cartan matrix $C=[c_{ij}]_{i,j=1}^{m}$ is defined by
\ben
c_{ij}=\begin{cases}
\phantom{2\tss}(\al_i,\al_j)
\qquad&\text{if}\quad i<m,\\[0.2em]
2\tss(\al_i,\al_j)
\qquad&\text{if}\quad i=m.
\end{cases}
\een
We will use notation \eqref{enise} -- \eqref{efexp}
along with
\ben
e^\circ_i(u)=\sum_{r=2}^{\infty}e_i^{(r)}u^{-r}\Fand
f^\circ_i(u)=\sum_{r=2}^{\infty}f_i^{(r)}u^{-r}.
\een

\bth\label{thm:dp}
The extended Yangian $\X(\osp_{1|2m})$ is generated by
the coefficients of the series
$h_i(u)$ with $i=1,\dots,m+1$, the series
$e_i(u)$, $f_i(u)$ with $i=1,\dots, m$, and the series $e_{mm'}(u)$, $f_{m'm}(u)$,
subject only to the following relations, where the indices
take all admissible values unless specified otherwise.
We have
\begin{align}
\label{hihj}
\big[h_i(u),h_j(v)\big]&=0, \\
\label{eifj}
\big[e_i(u),f_j(v)\big]&=\delta_{i\tss j}\ts\frac{k_i(u)-k_i(v)}{u-v}\ts (-1)^{\overline{i+1}}.
\end{align}
For all pairs $(i,j)$ except for $(m+1,m)$ we have
\begin{align}
\label{hiej}
\big[h_i(u),e_j(v)\big]&=-(\ve_i,\al_j)\ts
\frac{h_i(u)\tss\big(e_j(u)-e_j(v)\big)}{u-v},\\[0.4em]
\label{hifj}
\big[h_i(u),f_j(v)\big]&=(\ve_i,\al_j)\ts
\frac{\big(f_j(u)-f_j(v)\big)\tss h_i(u)}{u-v},
\end{align}
where $\ve_{m+1}:=0$, while
\begin{align}
\label{mohtej}
\big[h_{m+1}(u),e_m(v)\big]&
=h_{m+1}(u)\,\Big(\frac{e_m(u)-e_m(v)}{u-v}
-\frac{e_m(u-1/2)-e_m(v)}{u-v-1/2}\Big),\\[0.4em]
\label{mohtfj}
\big[h_{m+1}(u),f_m(v)\big]&=\Big({-}\frac{f_m(u) -f_m(v) }{u-v}
+\frac{f_m(u-1/2)-f_m(v)}{u-v-1/2}\Big)\,h_{m+1}(u).
\end{align}
For $i=1,\dots,m-1$ we have
\begin{align}
\label{eiei}
&\big[e_i(u),e_{i}(v)\big]=-
\frac{\big(e_{i}(u)-e_{i}(v)\big)^2}{u-v},\\[0.4em]
\label{fifi}
&\big[f_i(u),f_{i}(v)\big]=\frac{\big(f_{i}(u)-f_{i}(v)\big)^2}{u-v},
\end{align}
whereas
\begin{align}
\non
\big[e_m(u),e_m(v)\big]&=\frac{e_m(u)^2+e_{mm'}(u)-e_m(v)^2-e_{mm'}(v)}{u-v}\\
{}&+\frac{e_m(u)\tss e_m(v)-e_m(v)\tss e_m(u)}{2\tss(u-v)}
-\frac{\big(e_m(u)-e_m(v)\big)^2}{2\tss(u-v)^2},
\label{moeiei}
\end{align}
\begin{align}
\non
\big[f_m(u),f_m(v)\big]&=\frac{f_m(u)^2-f_{m'm}(u)-f_m(v)^2+f_{m'm}(v)}{u-v}\\
{}&-\frac{f_m(u)\tss f_m(v)-f_m(v)\tss f_m(u)}{2\tss(u-v)}
-\frac{\big(f_m(u)-f_m(v)\big)^2}{2\tss(u-v)^2}.
\label{mofifi}
\end{align}
For $i<j$ we have
\begin{align}
\label{eiej}
u\big[e^{\circ}_i(u),e_{j}(v)\big]-v\big[e_i(u),e^{\circ}_{j}(v)\big]
&=-(\al_i,\al_{j})\tss e_{i}(u)\tss e_{j}(v),\\[0.4em]
\label{fifj}
u\big[f^{\circ}_i(u),f_{j}(v)\big]-v\big[f_i(u),f^{\circ}_{j}(v)\big]
&=(\al_i,\al_{j})\tss f_{j}(v)\tss f_{i}(u).
\end{align}
Furthermore,
\begin{align}
\non
\big[e_m(u),e_{mm'}(v)\big]&=-\frac{\big(e_m(u)-e_m(v)\big)
\big(e_{mm'}(u)-e_{mm'}(v)\big)}{u-v}\\[0.4em]
{}&-\frac{e_m(u+1/2)-e_m(v)}{u-v+1/2}\ts e_m(u)^2
-\frac{e_{mm'}(u+1/2)-e_{mm'}(v)}{u-v+1/2}\ts e_m(u),
\label{moeieoo}
\end{align}
\begin{align}
\non
\big[f_m(u),f_{m'm}(v)\big]&=\frac{\big(f_{m'm}(u)-f_{m'm}(v)\big)
\big(f_m(u)-f_m(v)\big)}{u-v}\\[0.4em]
{}&-f_m(u)^2\ts\ts\frac{f_m(u+1/2)-f_m(v)}{u-v+1/2}
+f_m(u)\ts\frac{f_{m'm}(u+1/2)-f_{m'm}(v)}{u-v+1/2},
\non
\end{align}
and
\begin{align}
\big[e_{m-1}^{(1)},e_{mm'}(v)\big]&=e_{m-1}(v)\tss e_{mm'}(v)
+e_{mm'}(v)\ts e_{m-1}(v-3/2)\non\\[0.4em]
{}&+e_m(v)\ts \big[e_m^{(1)},e_{m-1}(v-3/2)\big]
-\frac12\ts\big[e_{mm'}^{(1)},e_{m-1}(v)+e_{m-1}(v-3/2)\big],
\label{emne}
\end{align}
\begin{align}
\big[f_{m-1}^{(1)},f_{m'm}(v)\big]&=-f_{m'm}(v)\tss f_{m-1}(v)-f_{m-1}(v-3/2)\ts f_{m'm}(v)
\non\\[0.4em]
{}&-\big[f_m^{(1)},f_{m-1}(v-3/2)\big]\ts f_m(v)
-\frac12\ts\big[f_{m'm}^{(1)},f_{m-1}(v)+f_{m-1}(v-3/2)\big].
\non
\end{align}
Finally, we have
the Serre relations
\begin{align}
\label{eSerre}
\sum_{\si\in\Sym_k}\big[e_{i}(u_{\si(1)}),
\big[e_{i}(u_{\si(2)}),\dots,\big[e_{i}(u_{\si(k)}),e_{j}(v)\big]\dots\big]\big]&=0,\\
\label{fSerre}
\sum_{\si\in\Sym_k}\big[f_{i}(u_{\si(1)}),
\big[f_{i}(u_{\si(2)}),\dots,\big[f_{i}(u_{\si(k)}),f_{j}(v)\big]\dots\big]\big]&=0,
\end{align}
for $i\ne j$  with
$k=1+c_{ij}$.
\eth

\smallskip

\bpf
Relations~\eqref{hihj} were pointed out in Sec.~\ref{sec:gd} as consequences of
\eqref{ilm} and \eqref{cuhh}. By the Poincar\'e--Birkhoff--Witt theorem,
the Yangian $\Y(\gl_{0|m})$ can be regarded
as the subalgebra of $\X(\osp_{1|2m})$ generated by the coefficients of the series
$t_{ij}(u)$ with $1\leqslant i,j\leqslant m$. Hence, the relations involving
the Gaussian generators belonging to this subalgebra follow from \cite[Thm.~5.2]{bk:pp}.
Furthermore,
the relations involving the series $f_i(u)$ and $f_{m'm}(u)$ follow from
their counterparts involving $e_i(u)$ and $e_{mm'}(u)$
due to the symmetry
provided by the anti-automorphism $\tau$ defined in \eqref{tauanti}
which acts on the generators by \eqref{taue}. Relations \eqref{mohtej}, \eqref{moeiei}
and \eqref{moeieoo} follow from the respective relations of Theorem~\ref{thm:odp}
via the embedding theorem \cite[Thm~3.1]{m:dt}. Namely, the embedding
$\X(\osp_{1|2})\hra \X(\osp_{1|2m})$ constructed in {\em loc. cit.} is consistent
with the Gauss decompositions of the generator matrices and for the images
of the Gaussian generators of $\X(\osp_{1|2})$ we have
\ben
h_1(u)\mapsto h_m(u),\quad h_2(u)\mapsto h_{m+1}(u),\quad e(u)\mapsto e_{m}(u),\quad
e_{11'}(u)\mapsto e_{mm'}(u);
\een
see \cite[Prop~4.2]{m:dt}. Similarly, by using
the embedding $\X(\osp_{1|4})\hra \X(\osp_{1|2m})$
(with $m\geqslant 2$) we derive \eqref{emne} from the first relation
in Corollary~\ref{cor:reid}, where we take into account \eqref{paret}
and the relation $e_{13}(v)=[e_{23}^{(1)},e_{12}(v)]$ in $\X(\osp_{1|4})$.

The Serre relations for the series $e_i(u)$ and $f_i(u)$ are implied by
the Serre relations in the Lie superalgebra $\osp_{1|2m}$
(see e.g.~\cite[Sec.~2.44]{fss:dl})
via the embedding \eqref{emb}. This follows by
the argument originated in the work
of Levendorski\u\i~\cite[Lem.~1.4]{l:gd} in the same way as for
the extended Yangian $\X(\osp_{N|2m})$ with $N\geqslant 3$; see \cite[Sec.~7]{m:dt}.
The remaining cases of \eqref{eifj}, \eqref{hiej} and \eqref{eiej}
are verified by applying the corresponding arguments used in the proof
of \cite[Thm~6.1]{m:dt}, which rely on Cor.~3.3 and Lem.~4.3 therein; cf.
\cite[Prop.~5.11 and 5.13]{jlm:ib}.

We thus have
a homomorphism
\beql{surjhom}
\wh \X(\osp_{1|2m})\to\X(\osp_{1|2m}),
\eeq
where $\wh \X(\osp_{1|2m})$ denotes the (abstract) algebra
with generators and relations as in the statement of the theorem
and the homomorphism
takes
the generators to the elements
of $\X(\osp_{1|2m})$ with the same name.
 We will show that
this homomorphism is surjective and injective.

To prove the surjectivity, note that by \eqref{defrel},
\beql{tijto}
\big[t_{ij}(u),t^{(1)}_{j\ts j+1}\big]=-t_{i\ts j+1}(u)
\eeq
for $1\leqslant i<j\leqslant m$, while
\beql{tjjp}
\big[t^{(1)}_{j\ts j+1},t_{i\ts (j+1)'}(u)\big]=-t_{i\tss j'}(u)
\eeq
for $1\leqslant i\leqslant j\leqslant m$. Relations \eqref{tijto}, \eqref{tjjp} and
their counterparts obtained by the application of the anti-automorphism \eqref{tauanti}
together with
the Poincar\'e--Birkhoff--Witt theorem for the extended Yangian $\X(\osp_{1|2m})$
imply that this algebra is generated by the coefficients of the series $t_{ij}(u)$
with $1\leqslant i,j\leqslant m+1$. Hence, due to the Gauss decomposition
\eqref{gd}, the algebra $\X(\osp_{1|2m})$ is generated
by the coefficients of the series $h_i(u)$ for $i=1,\dots,m+1$ together with
$e_{ij}(u)$ and $f_{ji}(u)$ for $1\leqslant i<j\leqslant m+1$.
Write \eqref{tijto} and \eqref{tjjp} in terms of the Gaussian generators
(cf. \cite[Sec.~5]{jlm:ib}) to get
\beql{eijoop}
\big[e_{ij}(u),e^{(1)}_{j\ts j+1}\big]=-e_{i\ts j+1}(u)
\Fand
\big[e^{(1)}_{j\ts j+1},e_{i\ts (j+1)'}(u)\big]=-e_{i\tss j'}(u)
\eeq
for $1\leqslant i<j\leqslant m$, and
\beql{eiipp}
\big[e^{(1)}_{i\ts i+1},e_{i\ts (i+1)'}(u)\big]
=-e_{i\tss i'}(u)-e_{i\ts i+1}(u)\ts e_{i\ts (i+1)'}(u)
\eeq
for $i=1,\dots,m$. These relations together with their counterparts
for the coefficients of the series $f_{ji}(u)$, which are obtained by applying
the anti-automorphism $\tau$ via \eqref{taue}, show that
the coefficients of the series $h_i(u)$ for $i=1,\dots,m+1$ and
$e_i(u)$, $f_i(u)$ for $i=1,\dots,m$ generate the algebra $\X(\osp_{1|2m})$
thus proving that the homomorphism \eqref{surjhom} is surjective.

Now we turn to proving
the injectivity of the homomorphism \eqref{surjhom}.
It was shown in the proof of \cite[Thm~6.1]{m:dt} that
the set of monomials in
the generators
$h_{i}^{(r)}$
with $i=1,\dots,m+1$ and $r\geqslant 1$,
and $e_{ij}^{(r)}$ and $f_{ji}^{(r)}$ with $r\geqslant 1$
and the conditions
\beql{condij}
i<j\leqslant i'\qquad\text{for}\quad i=1,\dots,m,
\eeq
taken in some fixed order
with the powers of odd generators
not exceeding $1$, is linearly independent in the extended Yangian
$\X(\osp_{1|2m})$.

Furthermore, working now in the algebra $\wh \X(\osp_{1|2m})$,
introduce its elements inductively, as the coefficients of the series $e_{ij}(u)$
for $i$ and $j$ satisfying
\eqref{condij} by setting $e_{i\ts i+1}^{(r)}=e_i^{(r)}$
for $i=1,\dots,m$ and using relations
\eqref{eijoop} and \eqref{eiipp}. The defining relations
show that the map
\beql{antit}
\tau:e_{i}(u)\mapsto f_{i}(u),\qquad f_{i}(u)\mapsto
-e_{i}(u)(-1)^{\overline{i+1}}\qquad\text{for}\quad i=1,\dots,m,
\eeq
and $\tau:h_i(u)\mapsto h_i(u)$ for $i=1,\dots,m+1$, defines an anti-automorphism of
the algebra $\wh \X(\osp_{1|2m})$. (We use the same symbol as in \eqref{taue}, but this
should not cause a confusion since it is used for a differently defined algebra.)
Apply this map to the relations defining
$e_{ij}(u)$ and use \eqref{taue} to get the definition of
the coefficients of the series $f_{ji}(u)$ subject to the same
conditions \eqref{condij}.
Since the images of the elements $h_{i}^{(r)}$, $e_{ij}^{(r)}$ and $f_{ji}^{(r)}$
of the algebra $\wh \X(\osp_{1|2m})$ under the homomorphism \eqref{surjhom}
coincide with the elements of the extended Yangian $\X(\osp_{1|2m})$ denoted
by the same symbols, the injectivity of the homomorphism \eqref{surjhom} will be proved by showing that
the algebra $\wh \X(\osp_{1|2m})$ is spanned by
monomials in these elements
taken in some fixed order.

Denote by $\wh \Ec$, $\wh \Fc$
and $\wh \Hc$ the subalgebras of $\wh \X(\osp_{1|2m})$ respectively
generated by all elements
of the form $e_{i}^{(r)}$, $f_{i}^{(r)}$ and $h_{i}^{(r)}$.
Define an ascending filtration
on $\wh \Ec$ by setting $\deg e_{i}^{(r)}=r-1$
and denote by $\gr\wh \Ec$ the corresponding associated graded algebra.
To establish the spanning property of
the monomials in the $e_{ij}^{(r)}$ in the subalgebra $\wh \Ec$, it will be enough to verify
the relations
\begin{multline}\label{bareijre}
\big[\eb_{i\tss j}^{\tss(r)},\eb_{k\tss l}^{\tss(s)}\big]=
\de_{k\tss j}\ts\eb_{i\tss l}^{\tss(r+s-1)}-\de_{i\tss l}\ts\eb_{kj}^{\tss(r+s-1)}\tss
(-1)^{(\bi+\bj)(\bk+\bl)}\\
-\de_{k\tss i'}\ts\eb_{j' l}^{\tss(r+s-1)}\tss (-1)^{\bi\tss\bj+\bi}\ts\ta_i\ta_j
+\de_{j'\tss l}\ts\eb_{k\tss i'}^{\tss(r+s-1)}\tss(-1)^{\bi\tss\bk+\bj\tss\bk+\bi+\bj}\ts\ta_i\ta_j,
\end{multline}
where $\eb_{ij}^{\tss(r)}$ denotes the image of the element $(-1)^{\bi}\ts e_{ij}^{(r)}$ in the
$(r-1)$-th component of $\gr\wh \Ec$ and we
extend the range of subscripts of
$\eb_{ij}^{\tss(r)}$ to all values $1\leqslant i<j\leqslant 1'$
by using the skew-symmetry conditions
\ben
\eb_{i\tss j}^{\tss(r)}=-\eb_{j'\tss i'}^{\tss(r)}\ts(-1)^{\bi\bj+\bi}\ts\ta_i\ta_j.
\een

First observe that
relations \eqref{bareijre} hold in the case $r=s=1$ because the defining relations of the theorem
restricted to the generators $e_i^{(1)}$ with $i=1,\dots,m$ reproduce the respective part of
the Serre--Chevalley presentation of the Lie superalgebra $\osp_{1|2m}$;
see e.g.~\cite[Sec.~2.44]{fss:dl}. Furthermore, the definitions
\eqref{eijoop} and \eqref{eiipp} of the elements $e_{ij}^{(r)}\in \wh \Ec$
imply the relations in the graded algebra $\gr\wh \Ec${\rm :}
\beql{greijoop}
\big[\eb^{\tss(r)}_{ij},\eb^{\tss(1)}_{j\ts j+1}\big]=\eb^{\tss(r)}_{i\ts j+1}\qquad\text{for}
\quad 1\leqslant i<j\leqslant m
\eeq
and
\beql{greiipp}
\big[\eb^{\tss(r)}_{i\ts (j+1)'},\eb^{\tss(1)}_{j\ts j+1}\big]=\eb^{\tss(r)}_{i\tss j'}
\ts(-1)^{\overline{j+1}}
\qquad\text{for}
\quad 1\leqslant i\leqslant j\leqslant m.
\eeq

Now write
\eqref{hiej} in terms of the coefficients by using \eqref{expafo} to get
\ben
\big[h_p^{(2)},e_j^{(r)}\big]=(\ve_p,\al_j)\ts\big(e_j^{(r+1)}
+h_p^{(1)}e_j^{(r)}\big).
\een
Extend the filtration on $\wh \Ec$ to the subalgebra $\wh\Bc$ of $\wh \X(\osp_{1|2m})$
generated by all elements
$e_{i}^{(r)}$ and $h_{i}^{(r)}$ by setting $\deg h_{i}^{(r)}=r-1$.
Hence, in the associated graded algebra $\gr\wh\Bc$ we have
\beql{hbasi}
\big[\hba_p^{\tss(2)},\eb_j^{\tss(r)}\big]=(\ve_p,\al_j)\ts \eb_j^{\tss(r+1)},
\eeq
where $\hba_p^{\tss(2)}$ is the image of $h_p^{(2)}$ in $\gr\wh\Bc$.

\ble\label{lem:heat}
For all $r,s\geqslant 1$ in the algebra $\gr\wh\Bc$ we
have
\beql{ebijj}
\big[\eb^{\tss(r)}_{ij},\eb^{\tss(s)}_{j\ts j+1}\big]=\eb^{\tss(r+s-1)}_{i\ts j+1}
\qquad\text{for}
\quad 1\leqslant i<j\leqslant m.
\eeq
Moreover, for all $p=1,\dots,m$ we also have
\beql{hbij}
\big[\hba_p^{\tss(2)},\eb_{i\tss j}^{\tss(r)}\big]
=(\ve_p,\al_{i\tss j})\ts \eb_{i\tss j}^{\tss(r+1)}
\qquad\text{for}
\quad 1\leqslant i<j\leqslant m+1.
\eeq
\ele

\bpf
Relation \eqref{eiej} implies
\beql{neie}
\big[\eb^{\tss(r+1)}_{j-1\ts j},\eb^{\tss(s)}_{j\ts j+1}\big]=
\big[\eb^{\tss(r)}_{j-1\ts j},\eb^{\tss(s+1)}_{j\ts j+1}\big]
\eeq
for all $r,s\geqslant 1$. This yields \eqref{ebijj} for $i=j-1$.
Continue by induction on $j-i$ (which is the length of the root $\al_{i\tss j}$)
and suppose that $j-i\geqslant 2$. Then by \eqref{greijoop},
\beql{comeeij}
\big[\eb^{\tss(r)}_{ij},\eb^{\tss(s)}_{j\ts j+1}\big]=
\big[[\eb^{\tss(r)}_{i\ts j-1},\eb^{\tss(1)}_{j-1\ts j}],\eb^{\tss(s)}_{j\ts j+1}\big].
\eeq
Observe that the commutator $[\eb^{\tss(r)}_{i\ts j-1},\eb^{\tss(s)}_{j\ts j+1}]$ is zero.
Indeed, by the first relation in \eqref{eijoop},
each element $e^{\tss(r)}_{i\ts j-1}\in\wh \Ec$ is a commutator of certain
coefficients of the series $e_{i}(u),\dots,e_{j-2}(u)$. However,
the commutator of each of these series with $e_j(u)$ is zero by
the Serre relations \eqref{eSerre}. Hence, using \eqref{neie},
we can write the commutator in \eqref{comeeij} as
\ben
\big[\eb^{\tss(r)}_{i\ts j-1},[\eb^{\tss(1)}_{j-1\ts j},\eb^{\tss(s)}_{j\ts j+1}]\big]
=\big[\eb^{\tss(r)}_{i\ts j-1},[\eb^{\tss(s)}_{j-1\ts j},\eb^{\tss(1)}_{j\ts j+1}]\big].
\een
Apply the Jacobi identity to this commutator.
By the induction hypothesis and \eqref{greijoop}, this equals
\ben
\big[\eb^{\tss(r+s-1)}_{i\tss j},\eb^{\tss(1)}_{j\ts j+1}\big]=\eb^{\tss(r+s-1)}_{i\ts j+1},
\een
as required, completing the proof of \eqref{ebijj}.

To verify \eqref{hbij}, use induction on $j-i$ with \eqref{hbasi} as the induction base.
For $j-i\geqslant 2$ write
\ben
\big[\hba_p^{\tss(2)},\eb_{i\tss j}^{\tss(r)}\big]=
\big[\hba_p^{\tss(2)},[\eb_{i\tss j-1}^{\tss(r)},\eb^{\tss(1)}_{j-1\ts j}]\big].
\een
By the induction hypothesis and \eqref{ebijj}, this equals
\begin{multline}
(\ve_p,\al_{i\tss j-1})\ts \big[\eb_{i\tss j-1}^{\tss(r+1)},\eb^{\tss(1)}_{j-1\ts j}\big]
+(\ve_p,\al_{j-1\tss j})\ts \big[\eb_{i\tss j-1}^{\tss(r)},\eb^{\tss(2)}_{j-1\ts j}\big]\\[0.4em]
{}=(\ve_p,\al_{i\tss j-1})\ts\eb_{i\tss j}^{\tss(r+1)}
+(\ve_p,\al_{j-1\tss j})\ts\eb_{i\tss j}^{\tss(r+1)}=(\ve_p,\al_{i\tss j})\ts\eb_{i\tss j}^{\tss(r+1)},
\non
\end{multline}
where we also used the root relation $\al_{i\tss j-1}+\al_{j-1\tss j}=\al_{i\tss j}$; see the notation
introduced in the beginning of Sec.~\ref{subsec:py}.
\epf

\ble\label{lem:heapr}
For all $r,s\geqslant 1$ and $1\leqslant i\leqslant j\leqslant m$ in the algebra $\gr\wh\Bc$ we
have
\beql{ebijjpr}
\big[\eb^{\tss(r)}_{i\ts (j+1)'},\eb^{\tss(s)}_{j\ts j+1}\big]=\eb^{\tss(r+s-1)}_{i\tss j'}
\ts(-1)^{\overline{j+1}}.
\eeq
Moreover, for all $p=1,\dots,m$ we also have
\beql{hbijpr}
\big[\hba_p^{\tss(2)},\eb^{\tss(r)}_{i\tss j'}\big]
=(\ve_p,\al_{i\tss j'})\ts \eb^{\tss(r+1)}_{i\tss j'}.
\eeq
\ele

\bpf
We will be proving both relations simultaneously by reverse induction on $j$
starting with $j=m$ (and then an inner induction on $i$).
In this case, relation \eqref{ebijjpr} with $i=m$ holds due to
\eqref{moeiei}, while using \eqref{hbasi} with $j=m$ we then derive \eqref{hbijpr}.
Now take $i=m-1$. Relation \eqref{emne} along with \eqref{ebijjpr} for $i=j=m$ give
\ben
\big[\eb^{\tss(1)}_{m-1\ts m},\eb^{\tss(s)}_{mm'}\big]=
-\big[\eb^{\tss(1)}_{m\tss m'},\eb^{\tss(s)}_{m-1\ts m}\big]
=\big[\eb^{\tss(s)}_{m-1\ts m},[\eb^{\tss(1)}_{m\ts m+1},\eb^{\tss(1)}_{m\ts m+1}]\big].
\een
Now take repeated commutators with $\hba_{m-1}^{\tss(2)}$ (which commutes
with $\eb^{\tss(s)}_{m\tss m'}$) to get
\ben
\big[\eb^{\tss(r)}_{m-1\ts m},\eb^{\tss(s)}_{mm'}\big]=
\big[\eb^{\tss(r+s-1)}_{m-1\ts m},[\eb^{\tss(1)}_{m\ts m+1},\eb^{\tss(1)}_{m\ts m+1}]\big]
=2\ts \eb^{\tss(r+s-1)}_{m-1\ts m'},
\een
where we also used \eqref{greijoop} and \eqref{greiipp}. Hence, by Lemma~\ref{lem:heat} the left hand side
of \eqref{ebijjpr} can be written as
\ben
\big[\eb^{\tss(r)}_{m-1\ts m+1},\eb^{\tss(s)}_{m\ts m+1}\big]=
\big[[\eb^{\tss(r)}_{m-1\ts m},\eb^{\tss(1)}_{m\ts m+1}],\eb^{\tss(s)}_{m\ts m+1}\big]
=-\big[\eb^{\tss(r+s-1)}_{m-1\ts m+1},\eb^{\tss(1)}_{m\ts m+1}\big]
+\big[\eb^{\tss(r)}_{m-1\ts m},\eb^{\tss(s)}_{mm'}\big],
\een
which coincides with $\eb^{\tss(r+s-1)}_{m-1\ts m'}$, as required. Relation \eqref{hbijpr}
in the case $i=m-1$ and $j=m$ follows by the same calculation as in the proof of
Lemma~\ref{lem:heat} with the use of the root relation
$\al_{m-1\ts m+1}+\al_{m\ts m+1}=\al_{m-1\ts m'}$.

Continue by reverse induction on $i$ and suppose that $i<m-1$. Invoking Lemma~\ref{lem:heat} again
and using the induction hypothesis, we get
\begin{multline}
\big[\eb^{\tss(r)}_{i\ts m+1},\eb^{\tss(s)}_{m\ts m+1}\big]=
\big[[\eb^{\tss(r)}_{i\ts i+1},\eb^{\tss(1)}_{i+1\ts m+1}],\eb^{\tss(s)}_{m\ts m+1}\big]
=\big[\eb^{\tss(r)}_{i\ts i+1},[\eb^{\tss(1)}_{i+1\ts m+1},\eb^{\tss(s)}_{m\ts m+1}]\big]\\[0.4em]
=\big[\eb^{\tss(r)}_{i\ts i+1},[\eb^{\tss(s)}_{i+1\ts m+1},\eb^{\tss(1)}_{m\ts m+1}]\big]
=\big[[\eb^{\tss(r)}_{i\ts i+1},\eb^{\tss(s)}_{i+1\ts m+1}],\eb^{\tss(1)}_{m\ts m+1}\big]
=\big[\eb^{\tss(r+s-1)}_{i\ts m+1},\eb^{\tss(1)}_{m\ts m+1}\big],
\non
\end{multline}
which equals $\eb^{\tss(r+s-1)}_{i\tss m'}$ by \eqref{greiipp}. This proves \eqref{ebijjpr}
in the case under consideration; relation \eqref{hbijpr} then also follows.

As a final step, continue by reverse induction on $j$ and suppose that
$1\leqslant i\leqslant j<m$. By \eqref{greiipp}
we have
\ben
\big[\eb^{\tss(r)}_{i\ts (j+1)'},\eb^{\tss(s)}_{j\ts j+1}\big]=
\big[[\eb^{\tss(r)}_{i\ts (j+2)'},\eb^{\tss(1)}_{j+1\ts j+2}],
\eb^{\tss(s)}_{j\ts j+1}\big]\ts(-1)^{\overline{j+2}}.
\een
Now observe that $[\eb^{\tss(r)}_{i\ts (j+2)'},\eb^{\tss(s)}_{j\ts j+1}]=0$. This relation
for $r=s=1$ holds as a particular case of \eqref{bareijre}. For arbitrary $r,s\geqslant 1$
the relation follows by taking repeated commutators with $\hba_p^{\tss(2)}$
for suitable values of $p$ by using \eqref{hbasi} and \eqref{hbijpr}; it suffices to take $p=i$ and
$p=j+1$. Hence by Lemma~\ref{lem:heat},
\begin{multline}
\big[\eb^{\tss(r)}_{i\ts (j+1)'},\eb^{\tss(s)}_{j\ts j+1}\big]=
\big[\eb^{\tss(r)}_{i\ts (j+2)'},[\eb^{\tss(1)}_{j+1\ts j+2},
\eb^{\tss(s)}_{j\ts j+1}]\big]\ts(-1)^{\overline{j+2}}=
\big[\eb^{\tss(r)}_{i\ts (j+2)'},[\eb^{\tss(s)}_{j+1\ts j+2},
\eb^{\tss(1)}_{j\ts j+1}]\big]\ts(-1)^{\overline{j+2}}\\[0.4em]
=\big[[\eb^{\tss(r)}_{i\ts (j+2)'},\eb^{\tss(s)}_{j+1\ts j+2}],
\eb^{\tss(1)}_{j\ts j+1}\big]\ts(-1)^{\overline{j+2}}
=\big[\eb^{\tss(r+s-1)}_{i\ts (j+1)'},\eb^{\tss(1)}_{j\ts j+1}\big]
=\eb^{\tss(r+s-1)}_{i\tss j'}
\ts(-1)^{\overline{j+1}},
\non
\end{multline}
where the last equality holds by \eqref{greiipp}, while the second last
equality is valid by the induction hypothesis.
This proves
\eqref{ebijjpr}, while \eqref{hbijpr} then follows by the same argument as in the proof
of Lemma~\ref{lem:heat}.
\epf

We will now complete the verification of \eqref{bareijre}.
Lemmas~\ref{lem:heat} and \ref{lem:heapr} imply the commutation relations
\ben
[\hba_p^{\tss(2)},\eb_{i\tss j}^{\tss(r)}]=(\ve_p,\al_{i\tss j})\ts \eb_{i\tss j}^{\tss(r+1)}
\een
for all positive roots $\al_{i\tss j}$.
Then the commutator of $\hba_p^{\tss(2)}$ with the left hand side of
\eqref{bareijre} equals
\ben
(\ve_p,\al_{i\tss j})\ts\big[\eb_{i\tss j}^{\tss(r+1)},\eb_{k\tss l}^{\tss(s)}\big]
+(\ve_p,\al_{k\tss l})\ts\big[\eb_{i\tss j}^{\tss(r)},\eb_{k\tss l}^{\tss(s+1)}\big].
\een
First consider fixed parameters $i<j\leqslant i'$ and $k<l\leqslant k'$
satisfying the following condition:
there exist two different values $p=a$ and $p=b$ such that
\ben
\begin{vmatrix}(\ve_a,\al_{i\tss j})&(\ve_a,\al_{k\tss l})\\
(\ve_b,\al_{i\tss j})&(\ve_b,\al_{k\tss l})
\end{vmatrix}\ne 0.
\een
In this case, starting with \eqref{bareijre} for $r=s=1$ and taking repeated
commutators of both sides with $\hba_a^{\tss(2)}$ and $\hba_b^{\tss(2)}$
we derive the required
relations for the super-commutators by solving the arising system of two
linear equations. For instance, starting from
$
[\eb_{i\tss j}^{\tss(r)},\eb_{i\tss j'}^{\tss(s)}]=
\eb_{i\tss i'}^{\tss(r+s-1)}
$
with $1\leqslant i<j\leqslant m$, we can take $p=i$ and $p=j$ to use the induction
step by solving the system
of equations
\ben
\bal
\big[\eb_{i\tss j}^{\tss(r+1)},\eb_{i\tss j'}^{\tss(s)}\big]
+\big[\eb_{i\tss j}^{\tss(r)},\eb_{i\tss j'}^{\tss(s+1)}\big]&=
2\tss\eb_{i\tss i'}^{\tss(r+s)},\\[0.4em]
\big[\eb_{i\tss j}^{\tss(r+1)},\eb_{i\tss j'}^{\tss(s)}\big]
-\big[\eb_{i\tss j}^{\tss(r)},\eb_{i\tss j'}^{\tss(s+1)}\big]&=
0.
\eal
\een

Consider now the remaining cases, where the above condition
on the determinant cannot be satisfied.
To verify that
$
[\eb_{i\tss j}^{\tss(r)},\eb_{i\tss j}^{\tss(s)}]=0
$
for $1\leqslant i<j\leqslant m$, note first that for $j=i+1$ this follows from
\eqref{eiei}. Furthermore, if $i<k<j$ for some $k$, then by the previously verified
cases of \eqref{bareijre}, we have
\ben
\big[\eb_{i\tss j}^{\tss(r)},\eb_{i\tss j}^{\tss(s)}\big]=
\big[\eb_{i\tss j}^{\tss(r)},[\eb_{i\tss k}^{\tss(s)},\eb_{k\tss j}^{\tss(1)}]\big]=0,
\een
as required. For the next case ($j=m+1$), observe that by \eqref{moeiei}
\ben
\big[\eb_{m\tss m+1}^{\tss(r)},\eb_{m\tss m+1}^{\tss(s)}\big]=\eb_{m\tss m'}^{\tss(r+s-1)}.
\een
Hence, for $i<m$ we have
\ben
\big[\eb_{i\tss m+1}^{\tss(r)},\eb_{i\tss m+1}^{\tss(s)}\big]
=\big[\eb_{i\tss m+1}^{\tss(r)},[\eb_{i\tss m}^{\tss(s)},\eb_{m\tss m+1}^{\tss(1)}]\big]
=\big[\eb_{i\tss m}^{\tss(s)},\eb_{m\tss i'}^{\tss(r)}\big]
=\eb_{i\tss i'}^{\tss(r+s-1)},
\een
thus verifying this case. Finally, for $1\leqslant i\leqslant j\leqslant m$ we have
\ben
\big[\eb_{i\tss j'}^{\tss(r)},\eb_{i\tss j'}^{\tss(s)}\big]=
\big[\eb_{i\tss j'}^{\tss(r)},[\eb_{i\tss m+1}^{\tss(s)},\eb_{m+1\tss j'}^{\tss(1)}]\big]=0,
\een
completing the verification of \eqref{bareijre}.

By applying
the anti-automorphism \eqref{antit}, we deduce
from the spanning property of the ordered monomials in the
elements $e_{ij}^{\tss(r)}$,
that the ordered monomials
in the elements $f_{ji}^{(r)}$
with the powers of odd generators
not exceeding $1$,
span the subalgebra $\wh\Fc$.
It is clear that the ordered monomials in $h_i^{(r)}$
span $\wh\Hc$. Furthermore,
by the defining relations of $\wh \X(\osp_{1|2m})$, the multiplication map
\ben
\wh\Fc\ot\wh\Hc\ot\wh \Ec\to \wh \X(\osp_{1|2m})
\een
is surjective. Therefore, ordering the elements
$h_{i}^{(r)}$, $e_{ij}^{(r)}$ and $f_{ji}^{(r)}$ in such a way that
the elements of $\wh\Fc$ precede the elements
of $\wh\Hc$, and the latter precede the elements of $\wh \Ec$,
we can conclude that the ordered monomials in these elements
with the powers of odd generators
not exceeding $1$,
span $\wh \X(\osp_{1|2m})$. This
proves that \eqref{surjhom} is an isomorphism.
\epf

\bre\label{rem:newarg}
The argument used for the verification of \eqref{bareijre} provides a different proof
of the respective relations in the (super) Yangians; cf. \cite{bk:pp}, \cite{g:gd} and \cite{jlm:ib}.
\qed
\ere

Let $\Ec$, $\Fc$
and $\Hc$ denote the subalgebras of $\X(\osp_{1|2m})$ respectively
generated by all elements
of the form $e_{i}^{(r)}$, $f_{i}^{(r)}$ and $h_{i}^{(r)}$.
Consider the generators
$h_{i}^{(r)}$
with $i=1,\dots,m+1$ and $r\geqslant 1$,
and $e_{ij}^{(r)}$ and $f_{ji}^{(r)}$ with $r\geqslant 1$
and conditions \eqref{condij}.
Suppose that the elements
$h_{i}^{(r)}$, $e_{ij}^{(r)}$ and $f_{ji}^{(r)}$ are ordered in such a way that
the elements of $\Fc$ precede the elements
of $\Hc$, and the latter precede the elements of $\Ec$.
The following is a version of the Poincar\'e--Birkhoff--Witt
theorem for the orthosymplectic Yangian.

\bco\label{cor:pbwdp}
The set of all ordered monomials in the elements
$h_{i}^{(r)}$ with $i=1,\dots,m+1$, and the elements
$e_{ij}^{(r)}$ and $f_{ji}^{(r)}$ with $r\geqslant 1$
with the powers of odd elements
not exceeding $1$ and
satisfying conditions \eqref{condij},
forms a basis of the algebra $\X(\osp_{1|2m})$.
\qed
\eco

We will now apply Theorem~\ref{thm:dp} to
deduce a Drinfeld-type presentation for the Yangian $\Y(\osp_{1|2m})$.
By making use of the series \eqref{defkn}, introduce
the elements $\ka^{}_{i\tss r}$, $\xi_{i\tss r}^{\pm}$ and $\xi_{r}^{\pm}$
of the algebra $\X(\osp_{1|2m})$
as the coefficients of the series
\ben
\ka^{}_i(u)=1+\sum_{r=0}^{\infty}\ka^{}_{i\tss r}\ts u^{-r-1},
\qquad
\xi_i^{\pm}(u)=\sum_{r=0}^{\infty}\xi_{i\tss r}^{\pm}\ts u^{-r-1}
\Fand
\xi^{\pm}(u)=\sum_{r=0}^{\infty}\xi_{r}^{\pm}\ts u^{-r-1}
\een
by setting
\ben
\bal
\ka^{}_i(u)&=k_i\big(u-(m-i)/2\big),\\[0.3em]
\xi^{+}_i(u)&=f_i\big(u-(m-i)/2\big),\\[0.3em]
\xi^{-}_i(u)&=-e_i\big(u-(m-i)/2\big),
\eal
\een
for $i=1,\dots,m$, and
\ben
\xi^{+}(u)=f_{m'm}(u),\qquad
\xi^{-}(u)=-e_{mm'}(u).
\een
Since these series are fixed by all
automorphisms \eqref{muf}, their coefficients
belong to the subalgebra $\Y(\osp_{1|2m})$ of the extended Yangian $\X(\osp_{1|2m})$.
We will use the abbreviation $\{a,b\}=ab+ba$.

\bco\label{cor:modpy}
The Yangian $\Y(\osp_{1|2m})$ is generated by
the coefficients of the series
$\ka^{}_i(u)$, $\xi_i^{\pm}(u)$ for $i=1,\dots,m$, and $\xi^{\pm}(u)$
subject only to the following relations, where the indices
take all admissible values unless specified otherwise.
We have
\begin{align}
\label{kikj}
\big[\kappa_i(u),\kappa_j(v)\big]&=0,\\
\label{xpixmj}
\big[\xi_{i}^{+}(u),\xi_{j}^{-}(v)\big]&=-\de_{i\tss j}\ts\frac{\kappa_i(u)-\kappa_i(v)}{u-v}.
\end{align}
For all pairs $(i,j)$ except for $(m,m)$ we have
\begin{align}
\label{kixpj}
\big[\kappa_i(u),\xi^{\pm}_j(v)\big]&={}\mp\frac{1}{2}\ts(\al_i,\al_j)\ts
\frac{\big\{\kappa_i(u),\xi^{\pm}_j(u)-\xi^{\pm}_j(v)\big\}}{u-v},\\
\label{xpixpj}
\big[\xi^{\pm}_i(u),\xi^{\pm}_{j}(v)\big]
+\tss\big[\xi^{\pm}_j(u),\xi^{\pm}_{i}(v)\big]&=
{}\mp\frac{1}{2}\ts(\al_i,\al_j)
\frac{\big\{\xi^{\pm}_i(u)-\xi^{\pm}_i(v),
\xi^{\pm}_j(u)-\xi^{\pm}_j(v)\big\}}{u-v},
\end{align}
while
\begin{align}\label{mkufv}
\big[\ka_m(u),\xi_{m}^{+}(v)\big]&=\Big(\frac{\xi_{m}^{+}(u-1/2) -\xi_{m}^{+}(v) }{3\tss(u-v-1/2)}
+\frac{2\tss\big(\xi_{m}^{+}(u+1)-\xi_{m}^{+}(v)\big)}{3\tss(u-v+1)}\Big)\ts \ka_m(u),\\[0.4em]
\label{mkuev}
\big[\ka_m(u),\xi_{m}^{-}(v)\big]
&=\ka_m(u)\ts\Big({-}\frac{\xi_{m}^{-}(u-1/2) -\xi_{m}^{-}(v) }{3\tss(u-v-1/2)}
-\frac{2\tss\big(\xi_{m}^{-}(u+1)-\xi_{m}^{-}(v)\big)}{3\tss(u-v+1)}\Big),
\end{align}
and
\begin{align}
\non
\big[\xi_{m}^{\pm}(u),\xi_{m}^{\pm}(v)\big]&=
\frac{\xi_{m}^{\pm}(u)^2-\xi^{\pm}(u)-\xi_{m}^{\pm}(v)^2+\xi_{m}^{\pm}(v)}{u-v}\\
{}&\mp\frac{\xi_{m}^{\pm}(u)\tss \xi_{m}^{\pm}(v)-\xi_{m}^{\pm}(v)\tss \xi_{m}^{\pm}(u)}{2\tss(u-v)}
-\frac{\big(\xi_{m}^{\pm}(u)-\xi_{m}^{\pm}(v)\big)^2}{2\tss(u-v)^2}.
\label{mmoeiei}
\end{align}
Furthermore,
\begin{align}
\non
\big[\xi_{m}^{+}(u),\xi^{+}(v)\big]&=\frac{\big(\xi^{+}(u)-\xi^{+}(v)\big)
\big(\xi_{m}^{+}(u)-\xi_{m}^{+}(v)\big)}{u-v}\\[0.4em]
{}&-\xi_{m}^{+}(u)^2\ts\ts\frac{\xi_{m}^{+}(u+1/2)-\xi_{m}^{+}(v)}{u-v+1/2}
+\xi_{m}^{+}(u)\ts\frac{\xi^{+}(u+1/2)-\xi^{+}(v)}{u-v+1/2},
\label{mmofifoo}
\end{align}
\medskip
\begin{align}
\non
\big[\xi_m^{-}(u),\xi^{-}(v)\big]&=-\frac{\big(\xi_m^{-}(u)-\xi_m^{-}(v)\big)
\big(\xi^{-}(u)-\xi^{-}(v)\big)}{u-v}\\[0.4em]
{}&+\frac{\xi_m^{-}(u+1/2)-\xi_m^{-}(v)}{u-v+1/2}\ts \xi_m^{-}(u)^2
-\frac{\xi^{-}(u+1/2)-\xi^{-}(v)}{u-v+1/2}\ts \xi_m^{-}(u)
\label{mmoeieoo}
\end{align}
and
\begin{align}
\big[\xi^+_{m-1, 0},\xi^+(v)\big]&=-\xi^+(v)\tss \xi^+_{m-1}(v+1/2)\tss
-\big[(\xi^+_{m\ts 0})^2,\xi^+_{m-1}(v+1/2)+\xi^+_{m-1}(v-1)\big]
\non\\[0.4em]
{}&-\big[\xi^+_{m\tss 0},\xi^+_{m-1}(v-1)\big]\ts \xi^+_m(v)-\xi^+_{m-1}(v-1)\ts \xi^+(v),
\label{xifmne}
\end{align}
\begin{align}
\big[\xi^-_{m-1, 0},\xi^-(v)\big]&=\xi^-_{m-1}(v+1/2)\tss\xi^-(v)
-\big[(\xi^-_{m\ts 0})^2,\xi^-_{m-1}(v+1/2)+\xi^-_{m-1}(v-1)\big]
\non\\[0.4em]
{}&-\xi^-_m(v)\ts\big[\xi^-_{m\tss 0},\xi^-_{m-1}(v-1)\big]+\xi^-(v)\ts\xi^-_{m-1}(v-1).
\label{xiemne}
\end{align}
Finally,
the Serre relations
\beql{Serrexipm}
\sum_{\si\in\Sym_k}\big[\xi_i^{\pm}(u_{\si(1)}),
\big[\xi_i^{\pm}(u_{\si(2)}),\dots,
\big[\xi_i^{\pm}(u_{\si(k)}),\xi_j^{\pm}(v)\big]\dots\big]\big]=0,
\eeq
hold for all $i\ne j$, where we set $k=1+c_{ij}$.
\eco

\bpf
The relations
are deduced from Theorem~\ref{thm:dp} by the arguments similar to those in
\cite{bk:pp}; see also \cite[Sec.~3.1]{m:yc}.
In particular, \eqref{kixpj} and \eqref{xpixpj} are essentially the Yangian relations of type $A$,
while \eqref{mkufv} -- \eqref{mmoeieoo} follow from
Corollary~\ref{cor:odpy}
via the embedding theorem. To illustrate, we will derive \eqref{xpixpj} with $j=i+1$
for $\xi^{-}_i(u)$
from the corresponding case of \eqref{eiej}. We can write the latter in the form
\beql{eieii}
(u-v)\ts\big[e_i(u),e_{i+1}(v)\big]
=-e_{i}(u)\tss e_{i+1}(v)+A(u)+B(v)
\eeq
for certain series $A(u)$ and $B(u)$ in $u^{-1}$. By setting $v=u+1/2$ we derive
\beql{aubu}
A(u)+B(u+1/2)=\frac12\ts\big\{e_i(u),e_{i+1}(u+1/2)\big\}.
\eeq
Writing \eqref{xpixpj} in terms of the series $e_i(u)$ and shifting the variables by
$u\mapsto u+(m-i)/2$ and $v\mapsto v+(m-i)/2$ we come to verifying the relation
\begin{multline}
\big[e_i(u),e_{i+1}(v+1/2)\big]
-\big[e_i(v),e_{i+1}(u+1/2)\big]\\[0.4em]
=\frac{1}{2}\ts
\frac{\big\{e_i(u)-e_i(v),
e_{i+1}(u+1/2)-e_{i+1}(v+1/2)\big\}}{u-v}.
\label{veryf}
\end{multline}
Multiply both sides by $u-v$ and write
\begin{multline}
(u-v)\ts \big[e_i(u),e_{i+1}(v+1/2)\big]\\[0.4em]
=(u-v-\frac12)\ts \big[e_i(u),e_{i+1}(v+1/2)\big]
+\frac12\ts \big[e_i(u),e_{i+1}(v+1/2)\big]
\non
\end{multline}
to apply \eqref{eieii} to the first summand on the right hand side.
After expanding the commutators and anti-commutators in the resulting expression,
we conclude that it holds due to relation \eqref{aubu}.

The decomposition \eqref{tensordecom} and formula \eqref{cu}
imply that the coefficients of the series generate
the Yangian $\Y(\osp_{1|2m})$. The completeness of the relations
is verified by using the automorphisms of the form \eqref{muf}
on the abstract algebra with the presentation given in the statement of the corollary
as with the case $m=1$; see the proof of Corollary~\ref{cor:odpy}.
\epf

The relations of Corollary~\ref{cor:modpy} can be written explicitly in terms of the
generators $\ka^{}_{i\tss r}$, $\xi_{i\tss r}^{\pm}$ and $\xi_{r}^{\pm}$
by using the expansion \eqref{expafo}.
Most of them have the same form as for the Yangian $\Y(\osp_{N|2m})$
with $N\geqslant 3$ (see \cite[Main Theorem]{m:dt}), but those involving
shifts in $u$ are more complicated because they require further expansions
of series of the form $(u+a)^{-r}$.

\subsection{Highest weight representations}
\label{subsec:re}

We will conclude with an application of the results of \cite{m:ry} and
give a description
of the finite-dimensional irreducible representations of the Yangian
$\Y(\osp_{1|2m})$ in terms of the presentation of Corollary~\ref{cor:odpy}.

A representation $L$ of the Yangian $\Y(\osp_{1|2m})$
is called a {\em highest weight representation}
if there
exists a nonzero vector
$\ze\in L$ such that $L$ is generated by $\ze$
and the following relations hold:
\beql{hwr}
\xi^+_{i}(u)\ts\ze=0
\Fand
\ka_{i}(u)\ts\ze=\mu_i(u)\ts\ze \qquad
\text{for} \quad i=1,\dots,m,
\eeq
for some formal series
\beql{mui}
\mu_i(u)=1+\mu_i^{(1)}u^{-1}+\mu_i^{(2)}u^{-2}+\dots,\qquad
\mu_i^{(r)}\in\CC.
\eeq
The vector $\ze$ is called the {\em highest vector}
of $L$, and the $m$-tuple
$\mu(u)=(\mu_1(u),\dots,\mu_m(u))$
is the {\em highest weight}
of $L$.

Given an arbitrary tuple $\mu(u)=(\mu_1(u),\dots,\mu_m(u))$ of formal series
of the form \eqref{mui},
the {\em Verma module} $M(\mu(u))$ is defined as the quotient of the algebra $\Y(\osp_{1|2m})$ by
the left ideal generated by all coefficients of the series $\xi^+_{i}(u)$
and $\ka_i(u)-\mu_i(u)$
for $i=1,\dots,m$. We will denote by $L(\mu(u))$ its irreducible quotient.
This is a highest weight representation with the highest weight $\mu(u)$.
The isomorphism
class of $L(\mu(u))$ is determined by $\mu(u)$. The following description is
analogous to the classification theorem of \cite{d:nr}.

\bpr\label{prop:fdhw}
Every finite-dimensional irreducible representation of the algebra $\Y(\osp_{1|2m})$
is isomorphic to $L(\mu(u))$ for a certain tuple
$\mu(u)$. Moreover,
the representation $L(\mu(u))$ of $\Y(\osp_{1|2m})$
is finite-dimensional if and only if there exist monic polynomials
$Q_1(u),\dots,Q_{m}(u)$ in the variable
$u$ such that
\beql{fdco}
\mu_i(u)=\frac{Q_i(u+1)}{Q_i(u)}\qquad\text{for}\quad i=1,\dots,m.
\eeq
All $m$-tuples of monic polynomials $\big(Q_1(u),\dots,Q_{m}(u)\big)$
arise in this way.
\epr

\bpf
All parts of the proposition follow from the main theorem of \cite{m:ry} via the isomorphism
between the presentations of the algebra $\Y(\osp_{1|2m})$ constructed
in the proofs of Theorem~\ref{thm:dp} and Corollary~\ref{cor:modpy}.
In the same way as for the Yangians associated with the classical
Lie algebras (see \cite{jlm:rq} for more details), one only needs to twist the Yangian action
to relate the parameters of the highest and lowest weight
representations. More precisely, the action of the extended Yangian
in terms of the $RTT$ presentation as defined in Sec.~\ref{sec:def},
should be twisted
by the automorphism
\ben
t_{ij}(u)\mapsto t_{ji}(-u)(-1)^{\bi\tss\bj+\bj}
\een
for the first condition in \eqref{hwr} to correspond to the highest weight condition of \cite{m:ry}.
\epf

\subsection*{Data Availability Statement}

All data is available within the article.

\subsection*{Compliance with Ethical Standards}
This work was supported by the Australian Research Council, grant DP180101825.
The authors have no competing interests to declare that are relevant to the content of this article.

\bigskip
\bigskip

\small
\noindent
School of Mathematics and Statistics\newline
University of Sydney,
NSW 2006, Australia\newline
alexander.molev@sydney.edu.au

\vspace{5 mm}

\noindent
Laboratoire de Physique Th\'{e}orique LAPTh,
CNRS and Universit\'{e} de Savoie\newline
BP 110, 74941 Annecy-le-Vieux Cedex, France\newline
eric.ragoucy@lapth.cnrs.fr

\end{document}